\documentclass[aos,MSNbibl,nameyear,dvips]{arximspdf}
\usepackage{algorithm}
\usepackage{graphicx}

\doi{10.1214/11-AOS886}
\volume{39}
\issue{3}
\pubyear{2011}
\firstpage{1776}
\lastpage{1802}

\makeatletter
\renewcommand{\epsilon}{\varepsilon}
\newcommand\E{E}
\newcommand{\Var}{\operatorname{Var}}

\newcommand\loglik{\ell}
\newcommand\VarMC{\Var_{\mathit{MC}}}
\newcommand\EMC{\E_{\mathit{MC}}}
\newcommand\genX{\tilde{X}}
\newcommand\genf{\tilde{f}}
\newcommand\genY{\tilde{Y}}
\newcommand\gendata{{y}^{*}}
\newcommand\gentheta{\tilde{\theta}}
\newcommand\genRX{\tilde{\RX}}
\newcommand\genRY{{\RY}}
\newcommand\barX{\overline{X}}
\newcommand\barY{\overline{Y}}
\newcommand\barf{\overline{f}}
\newcommand\barloglik{\overline{\loglik}}
\newcommand\compact{K}
\newcommand\testFunction{\phi}
\newcommand\chTheta{\breve{\Theta}}
\newcommand\chX{\breve{X}}
\newcommand\chY{\breve{Y}}
\newcommand\chSigma{\breve{\Sigma}}
\newcommand\chE{{\breve{\mathbb E}}}
\newcommand\chVar{\breve{\operatorname{Var}}}
\newcommand\Egg{{\breve{\mathbb E}_{\theta,\tau}}}
\newcommand\YAphi{\psi}
\newcommand\varc{{\overline{c}}}
\newcommand\RX{{\mathfrak{X}}}
\newcommand\RY{{\mathfrak{Y}}}
\newcommand\data{y^{*}}
\newcommand\Z{\zeta}
\newcommand\kA{{K_{1}}}
\newcommand\kB{{K_{3}}}
\newcommand\kC{{K_{2}}}
\newcommand\kD{{K_{4}}}
\newcommand\kE{{K_{6}}}
\newcommand\kF{{K_{7}}}
\newcommand\kH{{K_{5}}}
\newcommand\kI{{K_{8}}}
\newcommand\rE{{R_{2}}}
\newcommand\rF{{R_{1}}}
\newcommand\rB{{R_{3}}}
\newcommand\rC{{R_{4}}}
\newcommand\cDD{{C_{7}}}
\newcommand\cEE{{C_{8}}}
\newcommand\cFF{{C_{27}}}
\newcommand\cGG{{C_{1}}}
\newcommand\cHH{{C_{2}}}
\newcommand\cII{{C_{3}}}
\newcommand\cJJ{{C_{5}}}
\newcommand\cKK{{C_{4}}}
\newcommand\cLL{{C_{9}}}
\newcommand\cMM{{C_{10}}}
\newcommand\cMMM{{C_{11}}}
\newcommand\cNN{{C_{12}}}
\newcommand\cOO{{C_{14}}}
\newcommand\cPP{{C_{0}}}
\newcommand\cQQ{{C_{13}}}
\newcommand\cRR{{C_{15}}}
\newcommand\cSS{{C_{16}}}
\newcommand\cTT{{C_{6}}}
\newcommand\cA{{C_{26}}}
\newcommand\cB{{C_{21}}}
\newcommand\cD{{C_{17}}}
\newcommand\cS{{C_{25}}}
\newcommand\cT{{C_{23}}}
\newcommand\cU{{C_{24}}}
\newcommand\cV{{C_{22}}}
\newcommand\cW{{C_{18}}}
\newcommand\cY{{C_{19}}}
\newcommand\cZ{{C_{20}}}

\def\filter{{F}}
\def\predict{{P}}
\newcommand{\filttheta}{\breve{\theta}^F }
\newcommand{\predV}{\breve{V}^P }
\newcommand{\given}{\mid}
\newcommand{\kernel}{\kappa}
\newcommand{\equals}{=}
\newcommand{\prob}{{\mathbb{P}}}
\newcommand{\dimtheta}{p}
\newcommand{\R}{\mathbb{R}}
\newcommand{\ParScale}{\Sigma}
\newcommand{\transpose}{{\prime}}

\newtheorem{theorem}{Theorem}
\newtheorem{corollary}[theorem]{Corollary}
\newcommand{\fracd}[2]{({#1}/{#2})}

\setattribute{abstract}{width}{280pt}

\makeatother

\begin{document}
\begin{frontmatter}

\title{Iterated filtering\thanksref{TT1}}
\runtitle{Iterated filtering}

\begin{aug}
\author[A]{\fnms{Edward L.} \snm{Ionides}\corref{}\ead[label=e1]{ionides@umich.edu}},
\author[B]{\fnms{Anindya} \snm{Bhadra}\ead[label=e2]{tatar@umich.edu}},
\author[B]{\fnms{Yves}~\snm{Atchad\'{e}}\ead[label=e3]{yvesa@umich.edu}}
\and
\author[C]{\fnms{Aaron}~\snm{King}\ead[label=e4]{kingaa@umich.edu}}
\runauthor{Ionides, Bhadra, Atchad\'{e} and King}
\affiliation{University of Michigan and Fogarty International Center,
National Institutes of Health,
University of Michigan,
University of Michigan,
and~University of Michigan and Fogarty International Center, National~Institutes of Health}
\address[A]{E. L. Ionides\\
Department of Statistics\\
 University of Michigan\\
Ann Arbor, Michigan 48109\\
USA\\
and\\
Fogarty International Center\\
National Institutes of Health\\
Bethesda, Maryland 20892\\
 USA\\
\printead{e1}} 
\address[B]{A. Bhadra\\
Y. Atchad\'{e}\\
Department of Statistics\\
 University of Michigan\\
Ann Arbor, Michigan\\
USA\\
\printead{e2}\\
\hphantom{\textsc{E-mail:}} \printead*{e3}}
\address[C]{A. King\\
Department of Ecology and\\ \quad  Evolutionary Biology\\
  University of Michigan\\
Ann Arbor, Michigan\\
USA\\
and\\
Fogarty International Center\\
National Institutes of Health\\
Bethesda, Maryland 20892\\
 USA\\
\printead{e4}}
\end{aug}
\thankstext{TT1}{Supported by  NSF Grants DMS-08-05533 and EF-04-30120, the Graham
Environmental Sustainability Institute, the RAPIDD program of the
Science \& Technology Directorate, Department of Homeland Security, and
the Fogarty International Center, National Institutes of Health.}

\thankstext[]{}{This work was conducted as part of the Inference for Mechanistic Models Working Group
supported by the National Center for Ecological Analysis and Synthesis, a Center funded by NSF
Grant DEB-0553768, the University of California, Santa Barbara and the State of
California.}

\received{\smonth{9} \syear{2010}}
\revised{\smonth{3} \syear{2011}}

%
\begin{abstract}
Inference for partially observed Markov process models has been a
longstanding methodological challenge with many scientific and
engineering applications. Iterated filtering algorithms maximize the
likelihood function for partially observed Markov process models by
solving a recursive sequence of filtering problems. We present new
theoretical results pertaining to the convergence of iterated filtering
algorithms implemented via sequential Monte Carlo filters. This theory
complements the growing body of empirical evidence that iterated
filtering algorithms provide an effective inference strategy for
scientific models of nonlinear dynamic systems. The first step in our
theory involves studying a new recursive approach for maximizing the
likelihood function of a latent variable model, when this likelihood is
evaluated via importance sampling. This leads to the consideration of
an iterated importance sampling algorithm which serves as a simple
special case of iterated filtering, and may have applicability in its
own right.
\end{abstract}

\setattribute{keyword}{AMS}{MSC2010 subject classification.}

\end{frontmatter}

\section{Introduction} 

Partially observed Markov process (POMP) models are of widespread
importance throughout science and engineering.
As such, they have been studied\vadjust{\goodbreak} under various names including \textit
{state space models}
[\citet{durbin01}], \textit{dynamic models} [\citet
{west97}] and \textit{hidden Markov models} [\citet{cappe05}].
Applications include ecology [\citet{newman08}], economics
[\citet{fernandez07}], epidemiology [\citet{king08}],
finance [\citet{johannes08}], meteorology [\citet
{anderson07}], neuroscience [\citet{ergun07}] and target tracking
[\citet{godsill07}].

This article investigates convergence of a Monte Carlo technique for
estimating unknown parameters of POMPs, called \textit{iterated
filtering}, which was proposed by \citet{ionides06-pnas}.
Iterated filtering algorithms repeatedly carry out a filtering
procedure to explore the likelihood surface at increasingly local
scales in search of a maximum of the likelihood function.
In several case-studies, iterated filtering algorithms have been shown
capable of addressing scientific challenges in the study of infectious
disease transmission, by making likelihood-based inference
computationally feasible in situations where this was previously not
the case [\citet{king08};
\citet{breto09};
\citet{he10};
\citet{laneri10}].
The partially observed nonlinear stochastic systems arising in the
study of disease transmission and related ecological systems are a
challenging environment to test statistical methodology [\citet
{bjornstad01}], and many statistical methodologies have been tested on
these systems in the past fifty years
[e.g., \citet{cauchemez08-jrsi};
\citet{toni08};
\citet{keeling08};
\citet{ferrari08};
\citet{morton05};
\citet{grenfell02};
\citet{kendall99};
\citet{bartlett60};
\citet{bailey55}].
Since iterated filtering has already demonstrated potential for
generating state-of-the-art analyses on a major class of scientific
models, we are motivated to study its theoretical justification.
The previous theoretical investigation of iterated filtering, presented
by \citet{ionides06-pnas}, did not engage directly in the Monte Carlo
issues relating to practical implementation of the methodology.
It is relatively easy to check numerically that a maximum has been
attained, up to an acceptable level of Monte Carlo uncertainty, and
therefore one can view the theory of \citet{ionides06-pnas} as
motivation for an algorithm whose capabilities were proven by demonstration.
However, the complete framework presented in this article gives
additional insights into the potential capabilities, limitations and
generalizations of iterated filtering.

The foundation of our iterated filtering theory is a Taylor series
argument which we present first in the case of a general latent
variable model in Section~\ref{seciis}.
This leads us to propose and analyze a novel \textit{iterated importance
sampling} algorithm for maximizing the likelihood function of latent
variable models.
Our motivation is to demonstrate a relatively simple theoretical result
which is then extended to POMP models in Section~\ref{sectheory}.
However, this result also demonstrates the broader possibilities of the
underlying methodological approach.

The iterated filtering and iterated importance sampling algorithms that
we study have an attractive practical property that the model for the
unobserved process enters the algorithm only through the requirement
that realizations can be generated at arbitrary parameter values.
This property has been called \textit{plug-and-play} [\citet{breto09};
\citet{he10}] since it permits simulation code to be simply
plugged into the inference procedure.
A concept closely related to plug-and-play is that of \textit{implicit
models} for which the model is specified via an algorithm to generate
stochastic realizations [\citet{diggle84};
\citet{breto09}].
In particular, evaluation of the likelihood function for implicit
models is considered unavailable.
Implicit models arise when the model is represented by a ``black box''
computer program.
A scientist investigates such a model by inputting parameter values,
receiving as output from the ``black box'' independent draws from a
stochastic process, and comparing these draws to the data to make inferences.
For an implicit model, only plug-and-play statistical methodology can
be employed.
Other plug-and-play methods proposed for partially observed Markov
models include approximate Bayesian computations implemented via
sequential Monte Carlo [\citet{janeliu01};
\citet{toni08}], an asymptotically
exact Bayesian technique combining sequential Monte Carlo with Markov
chain Monte Carlo [\citet{andrieu10}], simulation-based forecasting
[\citet{kendall99}], and simulation-based spectral analysis
[\citet{reuman06}].
Further discussion of the plug-and-play property is included in the
discussion of Section~\ref{secdiscussion}.

\section{Iterated importance sampling} \label{seciis} 

Let $f^{}_{XY}(x,y; \theta)$ be the joint density of a pair of
random variables $(X,Y)$ depending on a parameter $\theta\in\R
^\dimtheta$.
We suppose that $(X,Y)$ takes values in some measurable space $\RX
\times\RY$, and $ f^{}_{XY}(x,y; \theta)$ is defined with respect
to some $\sigma$-finite product measure which we denote by $dx\,dy$.
We suppose that the observed data consist of a fixed value $\data\in
\RY$, with $X$ being unobserved.
Therefore, $\{f_{XY}(x,y; \theta), \theta\in\R^\dimtheta\}$
defines a general {latent variable} statistical model.
We write the marginal densities of $X$ and $Y$ as $f_X(x; \theta)$
and $f_Y(y; \theta)$, respectively.
The \textit{measurement model} is the conditional density of the observed
variable given the latent variable $X$, written as
$f_{Y|X}(y\given x; \theta)$.
The log likelihood function is defined as $\loglik(\theta)=\log
f_Y(\data; \theta)$.
We consider the problem of calculating the maximum likelihood estimate,
defined as $\hat\theta=\arg\max_\theta\loglik(\theta)$.

We consider an iterated importance sampling algorithm which gives a
plug-and-play approach to likelihood based inference for implicit
latent variable models, based on generating simulations at parameter
values in a neighborhood of the current parameter estimate to refine
this estimate.
This shares broadly similar goals with other Monte Carlo methods
proposed for latent variable models [e.g., \citet{johansen08};
\citet{qian06}], and in a more general context has
similarities with evolutionary optimization strategies [\citet{beyer01}].
We emphasize that the present motivation for proposing and studying
iterated importance sampling is to lay the groundwork for the results
on iterated filtering in Section~\ref{sectheory}.
However, the successes of iterated filtering methodology on POMP models
also raise the possibility that related techniques may be useful in
other latent variable situations.

We define the \textit{stochastically perturbed} model to be a triplet of
random variables $(\chX,\chY,\chTheta)$, with perturbation parameter
$\tau$ and parameter $\theta$, having a joint density on $\RX\times
\RY\times\R^\dimtheta$ specified as
%
%
\begin{equation}
g_{\chX,\chY,\chTheta}(x,y,\breve{\vartheta}; \theta,\tau
)=f^{}_{XY}(x,y; \breve{\vartheta})
\tau^{-\dimtheta}
\kernel_\tau
\bigl((\breve{\vartheta}-\theta) / \tau\bigr).
\label{eqgiis}
\end{equation}
Here, $\{\kernel_\tau, \tau>0\}$ is a collection of mean-zero
densities on $\R^\dimtheta$ (with respect to Lebesgue measure)
satisfying condition~\hyperlink{A:kappa:T=1}{(A1)} below:
\begin{longlist}
\item[(A1)]\hypertarget{A:kappa:T=1}
For each $\tau>0$, $\kernel_\tau$ is supported on a compact set
$K_0\subset\R^\dimtheta$ independent of $\tau$.
\end{longlist}
Condition~\hyperlink{A:kappa:T=1}{(A1)} can be satisfied by
the arbitrary selection of $\kernel_\tau$.
At first reading, one can imagine that $\kernel_\tau$ is fixed,
independent of $\tau$. However, the additional generality will be
required in Section~\ref{sectheory}.

We start by showing a relationship between conditional moments of
$\chTheta$ and the derivative of the log likelihood function, in
Theorem~\ref{thm1}.
We write $\Egg[\cdot]$ to denote expectation with respect to the
stochastically perturbed model.
We write $u\in\R^p$ to specify a column vector, with $u^\transpose$
being the transpose of $u$.
For a function $f=(f_1,\ldots,f_m)^\transpose\dvtx  \R^p\to\R^m$, we
write $\int f(u)\, du$ for the vector $ (\int f_1(u)\, du,\ldots,
\int f_m(u)\, du )^\transpose\in\R^m$;
For any function $f\dvtx \R^p\to\R$, we write $\nabla f(u)$ to denote the
column vector gradient of $f$, with $\nabla^2 f(u)$ being the second
derivative matrix.
We write $|\cdot|$ for the absolute value of a vector or the largest
absolute eigenvalue of a square matrix.
We write $B(r)=\{u\in\R^\dimtheta\dvtx  |u|\le r\}$ for the ball of
radius $r$ in $\R^\dimtheta$.
We assume the following regularity condition:

\begin{longlist}
\item[(A2)]\hypertarget{A:loglik:T=1}
$\loglik(\theta)$ is twice differentiable.
For any compact set $\kA\subset\R^\dimtheta$,
\[
\sup_{\theta\in\kA}|\nabla\loglik(\theta)|<\infty \quad
\mbox{and} \quad  \sup_{\theta\in\kA}|\nabla^2\loglik(\theta
)|<\infty.
\]
\end{longlist}

%
%
\begin{theorem}\label{thm1}
Assume conditions \hyperlink{A:kappa:T=1}{(\textup{A}1)}, \hyperlink
{A:loglik:T=1}{(\textup{A}2)}.
Let $h\dvtx  \R^\dimtheta\to\R^{m}$ be a measurable function possessing
constants $\alpha\geq0$, $c>0$ and $\epsilon>0$ such that, whenever
$u\in B(\epsilon)$,
%
%
\begin{equation}\label{condh}
|h(u)|\leq c |u|^\alpha.
\end{equation}
Define $\tau_0=\sup\{\tau\dvtx K_0\subset B(\epsilon/\tau)\}$.
For any compact set $\kC\subset\R^\dimtheta$ there exists $\cGG
<\infty$ and a positive constant $\tau_1\le\tau_0$ such that, for
all $0<\tau\le\tau_1$,
%
%
\begin{eqnarray}
\label{eq555}
&&\sup_{\theta\in\kC} \biggl|
\Egg[h(\chTheta-\theta)\given\chY\equals\data]-
\int
h(\tau u) \kernel_\tau(u)\,du\nonumber\\
&& \hspace*{76pt} {} - \tau\biggl\{ \int h(\tau u)
u^\transpose\kernel_\tau(u)\,du \biggr\} \nabla\loglik(\theta)
\biggr|
\\
&& \qquad \le\cGG\tau^{2+\alpha}.
\nonumber
\end{eqnarray}
\end{theorem}
\begin{pf}
Let $g_{\chY\vert\chTheta}(y \given\breve{\vartheta}; \theta
,\tau)$
denote the conditional density of $\chY$ given $\chTheta$.
We note that $g_{\chY\vert\chTheta}$ does not depend on either $\tau
$ or $\theta$, and so we omit\vspace*{-1pt} these dependencies below.
Then, $g_{\chY\vert\chTheta}(\data\given\breve{\vartheta})=\int
f_{XY}(x,\data; \breve{\vartheta})\, dx=\exp(\loglik(\breve
{\vartheta}) )$.
Changing\vspace*{-1pt} variable to $u=(\breve{\vartheta}-\theta) / \tau$, we calculate
%
%
\begin{eqnarray}\label{eq113}
\Egg[h(\chTheta-\theta)\given\chY\equals\data]
&=&
\frac{\int h(\breve{\vartheta}-\theta)g_{\chY|\chTheta
}(\data\given
\breve{\vartheta}) \tau^{-\dimtheta} \kernel_\tau((\breve
{\vartheta}-\theta)\big/ \tau)\, d\breve{\vartheta}}{\int g_{\chY
|\chTheta}(\data
\given\breve{\vartheta}) \tau^{-\dimtheta} \kernel_\tau((\breve
{\vartheta}-\theta)\big/ \tau)\, d\breve{\vartheta}}
\nonumber
\\[-8pt]
\\[-8pt]
&=&
\frac{\int h(\tau u) \exp\{\loglik(\theta+\tau u)-\loglik(\theta
)\} \kernel_\tau(u)\, du}{\int\exp\{\loglik(\theta+\tau
u)-\loglik(\theta)\} \kernel_\tau(u)\, du}.
\nonumber
\end{eqnarray}
Applying the Taylor expansion $e^x=1+x+ (\int
_0^1(1-t)e^{tx}\,dt )x^2$ to the numerator of (\ref{eq113}) gives
%
%
\begin{eqnarray} \label{eqtaylor1thm1}
&&\int\exp\{\ell(\theta+ \tau u)-\ell(\theta) \}h(\tau
u)\kappa_{\tau}(u)\,du\nonumber\\
&& \qquad =\int h(\tau u)\kappa_\tau(u)\,du + \int
\{\ell(\theta+\tau u)-\ell(\theta) \}h(\tau u)\kappa
_\tau(u)\,du\\
&& \qquad  \quad {}+\int\biggl(\int_0^1(1-t)e^{t(\ell(\theta+\tau u)-\ell(\theta
))}\,dt \biggr) \{\ell(\theta+\tau u)-\ell(\theta) \}^2
h(\tau u)\kappa_\tau(u)\, du.
\nonumber
\end{eqnarray}
We now expand the second term on the right-hand side of (\ref{eqtaylor1thm1}) by making use of the Taylor expansion
%
%
\begin{eqnarray}
\label{taylorloglik}
\ell(\theta+\tau u)-\ell(\theta)&=&\tau u^\transpose\int
_0^1\nabla\ell(\theta+t\tau u)\, dt\nonumber
\\[-8pt]
\\[-8pt] &=&\tau u^\transpose\nabla
\ell(\theta) +\tau u^\transpose\int_0^1 \{\nabla\ell
(\theta+t\tau u)-\nabla\ell(\theta) \}\,dt
\nonumber
\end{eqnarray}
and defining $\YAphi_{\tau,h}(\theta)=\int h(\tau u) \kappa_\tau
(u)\, du +\tau\{ \int h(\tau u)u^\transpose\kappa_\tau(u)\,
du \} \nabla\ell(\theta)$.
This allows us to rewrite (\ref{eqtaylor1thm1}) as
%
%
\begin{equation}
\label{eq1}
\int\exp\{ \ell(\theta+ \tau u)-\ell(\theta) \}
h(\tau u) \kappa_{\tau}(u)\, du
= \YAphi_{\tau,h}(\theta)+ \rF(\theta,\tau),
\end{equation}
where
\begin{eqnarray*}
&&\rF(\theta,\tau)\\
&& \qquad = \tau\int h(\tau u) u^\transpose\biggl(\int
_0^1 \{\nabla\ell(\theta+t\tau u)-\nabla\ell(\theta) \}\,dt \biggr)\kappa_\tau(u)\, du
\\
&& \qquad  \quad  {}+\int\biggl(\int_0^1(1-t) e^{t(\ell(\theta+\tau
u)-\ell(\theta))}\,dt \biggr) \{\ell(\theta+\tau u)-\ell(\theta
) \}^2h(\tau u) \kappa_\tau(u)\, du.
\end{eqnarray*}
%
As a consequence of \hyperlink{A:loglik:T=1}{(A2)}, we have
%
%
\begin{equation} \label{bound83}
\sup_{\theta\in\kC}\sup_{u\in K_0}\sup_{t\in[0,\tau_0]}
\bigl(|\nabla\ell(\theta+t u)|+|\nabla^2\ell(\theta+t u)|
\bigr)<\infty.
\end{equation}
Combining (\ref{bound83}) with the assumption that $K_0\subset
B(\epsilon/ \tau)$, we deduce the existence of a finite constant
$\cHH$ such that
%
%
\begin{equation}\label{eq1remainder}
\sup_{\theta\in\kC}|\rF(\theta,\tau)|\leq\cHH\tau
^{2+\alpha}.
\end{equation}
We bound the denominator of (\ref{eq113}) by considering the
special case of (\ref{eq1}) in which $h$ is taken to be the unit
function, $h(u)=1$.
Noting that $\int u\kappa_{\tau}(u)\,du=0$, we see that $\YAphi_{\tau
,1}(\theta)=1$ and so (\ref{eq1}) yields
\[
\int\exp\{\ell(\theta+ \tau u)-\ell(\theta) \}\kappa
_{\tau}(u)\,du
= 1+\rE(\theta,\tau)
\]
with $\rE(\theta,\tau)$ having a bound
%
%
\begin{equation}
\label{eq2remainder}
\sup_{\theta\in\kC}|\rE(\theta,\tau)|\leq\cII\tau^{2}
\end{equation}
for some finite constant $\cII$.
We now note the existence of a finite constant $\cKK$ such that
%
%
\begin{equation}\label{eq3}
\sup_{\theta\in\kC}|\YAphi_{\tau,h}(\theta)|\leq\cKK\tau
^\alpha
\end{equation}
implied by (\ref{condh}), \hyperlink{A:kappa:T=1}{(A1)} and
\hyperlink{A:loglik:T=1}{(A2)}.
Combining (\ref{eq1remainder}), (\ref{eq2remainder}) and (\ref{eq3}) with the identity
\[
\Egg\bigl(h(\chTheta-\theta)\given\chY\equals\data
\bigr)-\YAphi_{\tau,h}(\theta)
=
\frac{\rF(\theta,\tau)-\rE(\theta,\tau) \YAphi_{\tau
,h}(\theta)}{1+\rE(\theta,\tau)},
\]
and requiring $\tau_1<(2\cII)^{-1/2}$,
we obtain that
%
%
\begin{equation} \label{eqYAbound}
\sup_{\theta\in\kC} \bigl|\Egg\bigl(h(\chTheta-\theta)\given
\chY\equals\data\bigr)-\YAphi_{\tau,h}(\theta) \bigr| \leq
\cJJ\tau^{2+\alpha}
\end{equation}
for some finite constant $\cJJ$.
Substituting the definition of $\YAphi_{\tau,h}(\theta)$ into (\ref{eqYAbound}) gives (\ref{eq555}) and hence completes the proof.
\end{pf}

One natural choice is to take $h(u)=u$ in Theorem~\ref{thm1}.
Supposing that $\kernel_\tau$ has associated positive definite
covariance matrix $\chSigma$, independent of $\tau$, this leads to an
approximation to the derivative of the log likelihood given by
%
%
\begin{equation}
\label{eqabcde}\label{eq13}
\bigl|
\nabla\loglik(\theta) -
\{
(
\tau^2 \chSigma
)^{-1}
(
\Egg[\chTheta\given\chY\equals\data]-\theta
)
\}
\bigr|
\le\cTT\tau
\end{equation}
for some finite constant $\cTT$, with the bound being uniform over
$\theta$ in any compact subset of $\R^\dimtheta$.
The quantity $\Egg[\chTheta\given\chY\equals\data]$
does not usually have a closed form, but a plug-and-play Monte Carlo
estimate of it is available by importance sampling, supposing that one
can draw from $f_{X}(x; \theta)$ and evaluate $f_{Y|X}(\data
\given x; \theta)$.
Numerical approximation of moments is generally more convenient than
approximating derivatives, and this is the reason that the relationship
in (\ref{eqabcde}) may be useful in practice.
However, one might suspect that there is no ``free lunch'' and
therefore the numerical calculation of the left-hand side of (\ref{eqabcde}) should become fragile as $\tau$ becomes small.
We will see that this is indeed the case, but that iterated importance
sampling methods mitigate the difficulty to some extent by averaging
numerical error over subsequent iterations.

%
%
\begin{algorithm}
\def\qsp{\hspace{0.5cm}}
\textbf{Input:}
\begin{itemize}
\item Latent variable model
described by a latent variable density $f_{X}(x; \theta)$,
measurement model $f_{Y|X}(y\given x; \theta)$, and data
$\data$.
\item Perturbation density $\kernel$ having compact support, zero mean
and\vspace*{1pt} positive definite covariance matrix $\chSigma$.
%
\item Positive sequences $\{\tau_m\}$ and $\{a_m\}$
\item Integer sequence of Monte Carlo sample sizes, $\{J_m\}$
\item Initial parameter estimate, $\hat\theta_1$
\item Number of iterations, $M$
\end{itemize}

\textbf{Procedure:}
\begin{enumerate}[3]
\item[1] for $m$ in $1\dvtx M$
\item[2]\qsp for $j$ in $1\dvtx J_m$
\item[3]\qsp\qsp draw $Z_{j,m}\sim\kernel(\cdot)$ and set
$\chTheta_{j,m}=\hat\theta_{m}+\tau_mZ_{j,m}$
\item[4]\qsp\qsp draw $\chX_{j,m}\sim f_{X}(\cdot; \chTheta_{j,m})$
\item[5]\qsp\qsp set \rule{0mm}{4mm}$w_{j}=f_{Y|X}(\data
\given
\chX_{j,m}; \chTheta_{j,m})$
\item[6]\qsp end for
\item[7]\hypertarget{alg:iis:Dm}\qsp calculate $D_m = \tau
_m^{-2}\chSigma
^{-1} \{
(\sum_{j=1}^{J_m} w_{j} )^{-1}
(\sum_{j=1}^{J_m} w_{j}\chTheta_{j,m} ) - \hat\theta
_{m} \}$
\item[8]\qsp update estimate: $\hat\theta_{m+1}=\hat\theta_{m}+a_mD_m$
\item[9] end for
\end{enumerate}

\textbf{Output:}
\begin{itemize}
\item parameter estimate $\hat\theta_{M+1}$
\end{itemize}
\caption{A basic iterated importance sampling procedure.
The Monte Carlo random variables required at each iteration are
presumed to be drawn independently.
Theorem~\protect\ref{thsaiis} gives sufficient conditions for $\hat\theta
_M$ to converge to the maximum likelihood estimate as $M\to\infty$.}
\label{algiis}
\end{algorithm}

A trade-off between bias and variance is to be expected in any
Monte Carlo numerical derivative, a classic example being the
Kiefer--Wolfowitz algorithm [\citet{kiefer52};
\citet{spall03}].
Algorithms which are designed to balance such trade-offs have been
extensively studied under the label of \textit{stochastic approximation}
[\citet{kushner03};
\citet{spall03};
\citet{andrieu05}].
Algorithm~\ref{algiis} is an example of a basic stochastic
approximation algorithm taking advantage of (\ref{eqabcde}).
As an alternative, the derivative approximation in (\ref{eq13}) could be
combined with a stochastic line search algorithm.
In order to obtain the plug-and-play property, we consider an algorithm
that draws from $f_X(x; \theta)$ for iteratively selected values of
$\theta$.
This differs from other proposed iterative importance sampling
algorithms which aim to construct improved importance sampling
distributions [e.g., \citet{celeux06}].
In principle, a procedure similar to Algorithm~\ref{algiis} could
take advantage of alternative choices of importance sampling
distribution: the fundamental relationship in Theorem~\ref{thm1} is
separate from the numerical issues of computing the required
conditional expectation by importance sampling.
Theorem~\ref{thsaiis} gives sufficient conditions for the
convergence of Algorithm~\ref{algiis} to the maximum likelihood estimate.
To control the variance of the importance sampling weights, we suppose:
\begin{longlist}
\item[(A3)]\hypertarget{A:fbound:T=1}
For any compact set $\kB\subset\R^\dimtheta$,
\[
\sup_{\theta\in\kB, x\in\RX} f_{Y|X}(\data\given
x; \theta)<\infty.
\]
\end{longlist}
We also adopt standard sufficient conditions for stochastic
approximation methods:
\begin{longlist}[(B2)]
\item[(B1)]\hypertarget{C1}
Define $\Z(t)$ to be a solution to $d\Z/dt = \nabla\loglik(\Z(t))$.
Suppose that $\hat\theta$ is an \textit{asymptotically stable
equilibrium point}, meaning that
(i) for every $\eta>0$ there exists a $\delta(\eta)$ such that $|\Z
(t)-\hat\theta|\le\eta$ for all $t>0$ whenever $|\Z(0)-\hat\theta
|\le\delta$, and (ii) there exists a $\delta_0$ such that $\Z(t)\to
\hat\theta$ as $t\to\infty$ whenever $|\Z(0)-\hat\theta|\le
\delta_0$.
\item[(B2)]
\hypertarget{C2}
With probability one, $\sup_m |\hat\theta_m| <\infty$. Further,
$\hat\theta_m$ falls infinitely often into a compact subset of $\{\Z
(0)\dvtx \lim_{t\to\infty}\Z(t)=\hat\theta\}$.
\end{longlist}
 Conditions \hyperlink{C1}{(B1)} and \hyperlink{C2}{(B2)} are the
basis of the classic results of \citet{kushner78}.
Although research into stochastic approximation theory has continued
[e.g., \citet{kushner03};
\citet{andrieu05};
\citet{maryak08}],
\hyperlink{C1}{(B1)} and \hyperlink{C2}{(B2)} remain a textbook approach
[\citet{spall03}].
The relative simplicity and elegance of \citet{kushner78} makes an
appropriate foundation for investigating the links between iterated
filtering, sequential Monte Carlo and stochastic approximation theory.
There is, of course, scope for variations on our results based on the
diversity of available stochastic approximation theorems.
Although neither \hyperlink{C1}{(B1)} and \hyperlink{C2}{(B2)} nor
alternative sufficient conditions are easy to verify,
stochastic approximation methods have nevertheless been found effective
in many situations.
Condition \hyperlink{C2}{(B2)} is most readily satisfied if
$\hat\theta_m$ is constrained to a neighborhood in which $\hat\theta
$ is a unique local maximum, which gives a guarantee of local rather
than global convergence.
Global convergence results have been obtained for related stochastic
approximation procedures [\citet{maryak08}] but are beyond the
scope of
this current paper.
The rate assumptions in Theorem~\ref{thsaiis} are satisfied, for
example, by\vspace*{1pt} $a_m = m^{-1}$, $\tau_m^2 = m^{-1}$ and $J_m = m^{(\delta
+1/2)}$ for $\delta>0$.
%
%
\begin{theorem}\label{thsaiis}
Let $\{a_m\}$, $\{\tau_m\}$ and $\{J_m\}$ be positive sequences with
$\tau_m\to0$, $ J_m\tau_m \to\infty$, $a_m\to0$, $\sum_m
a_m=\infty$ and $\sum_m a_m^2J_m^{-1}\tau_m^{-2} <\infty$.
Let $\hat\theta_m$ be defined via Algorithm~\ref{algiis}.
Assuming \hyperlink{A:kappa:T=1}{(\textup{A}1)}--\hyperlink
{A:fbound:T=1}{(\textup{A}3)} and
\hyperlink{C1}{(\textup{B}1)} and \hyperlink{C2}{(\textup{B}2)},\break $\lim_{m\to\infty}
\hat\theta_m=\hat\theta$
with probability one.
\end{theorem}

\begin{pf}
The quantity $D_m$ in line~\hyperlink{alg:iis:Dm}{7} of Algorithm~\ref{algiis}
is a self-normalized Monte Carlo importance sampling estimate of
\[
(\tau_m^2 \chSigma)^{-1} \{{\mathbb E}_{\hat\theta
_{m},\tau_m} [\chTheta\given\chY\equals\data]-\hat
\theta_{m} \}.
\]
We can therefore apply Corollary~\ref{thcrisancor} from
Section~\ref{appsmc}, writing $\EMC$ and $\VarMC$ for the
Monte Carlo expectation and variance resulting from carrying out
Algorithm~\ref{algiis} conditional on the data $\data$.
This gives
%
%
\begin{eqnarray}
\label{isbias}
&&\bigl| \EMC\bigl[
D_m - (\tau_m^2 \chSigma)^{-1} \{{\mathbb E}_{\hat
\theta_{m},\tau_m} [\chTheta\given\chY\equals\data
]-\hat\theta_{m} \}
\bigr] \bigr|\nonumber
\\[-8pt]
\\[-8pt]
&& \qquad \le \frac{\cDD
(
\sup_{x\in\RX} f_{Y|X}(\data\given x; \hat\theta_{m})
)^2
}{
\tau_mJ_m
(
f_{Y}(\data; \hat\theta_{m})
)^2
},
\nonumber
\end{eqnarray}
\begin{equation}
 \label{isvar}
|
\VarMC(
D_m
)
|
\le \frac{\cEE(
\sup_{x\in\RX} f_{Y|X}(\data\given x; \hat\theta_{m})
)^2
}{
\tau_m^2J_m
(
f_{Y}(\data; \hat\theta_{m})
)^2
}
\end{equation}
for finite constants $\cDD$ and $\cEE$ which do not depend on $J$,
$\hat\theta_{m}$ or $\tau_m$.
Having assumed the conditions for Theorem~\ref{thm1}, we see from
(\ref{eqabcde}) and (\ref{isbias}) that $D_m$ provides an
asymptotically unbiased Monte Carlo estimate of $\nabla\loglik(\hat
\theta_{m})$ in the sense of \hyperlink{C5}{(B5)} of
Theorem~\ref{thspall} in Section~\ref{appsa}.
In addition, (\ref{isvar}) justifies \hyperlink{C4}{(B4)}
of Theorem~\ref{thspall}. The remaining conditions of Theorem~\ref{thspall} hold by hypothesis.
\end{pf}

\section{Iterated filtering for POMP models} 
\label{sectheory}

Let $\{X(t),t\in T\}$ be a Markov process [\citet{rogers94}] with $X(t)$
taking values in a measurable space~$\RX$.
The time index set $T\subset\R$ may be an interval\vadjust{\goodbreak} or a discrete set,
but we are primarily concerned with a finite subset of times
$t_1<t_2<\cdots<t_N$ at which $X(t)$ is observed, together with some
initial time $t_0<t_1$.
We write $X_{0\dvtx N}=(X_0,\ldots,X_N)= (X(t_0),\ldots,X(t_N) )$.
We write $Y_{1\dvtx N}=(Y_1,\ldots,Y_N)$ for a sequence of random variables
taking values in a measurable space $\RY^N$.
We assume that $X_{0\dvtx N}$ and $Y_{1\dvtx N}$ have a joint density
$f_{X_{0\dvtx N},Y_{1\dvtx N}}(x_{0\dvtx n},y_{1\dvtx n}; \theta)$ on $\RX^{N+1}\times
\RY^N$, with $\theta$ being an unknown parameter in $\R^{\dimtheta}$.
A POMP model may then be specified by an initial density
$f_{X_0}(x_0; \theta)$, conditional transition densities
$f_{X_n|X_{n-1}} (x_n\given x_{n-1} ; \theta)$ for $1\le n\le
N$, and the conditional densities of the observation process which are
assumed to have the form $f_{Y_n|Y_{1\dvtx n-1},X_{0\dvtx n}} (y_n\given
y_{1\dvtx n-1}, x_{0\dvtx n}; \theta) = f_{Y_n|X_n}(y_n\given x_n ;
\theta)$.
We use subscripts of $f$ to denote the required joint and conditional densities.
We write $f$ without subscripts to denote the full collection of
densities and conditional densities, and we call such an $f$ the
\textit{generic density} of a POMP model.
The data are a sequence of observations by $\data_{1\dvtx N}=(\data
_1,\ldots,\data_N)\in\RY^N$, considered as fixed.
We write the log likelihood function of the data for the POMP model as
$\loglik_N(\theta)$ where
\[
\loglik_n(\theta)=\log\int f_{X_0}(x_0; \theta)\prod_{k=1}^n
f_{X_k,Y_k |X_{k-1}}(x_k,\data_k\given x_{k-1} ; \theta) \, dx_{0\dvtx n}
\]
for $1\le n\le N$.
Our goal is to find the maximum likelihood estimate, $\hat\theta=\arg
\max_{\vartheta}\loglik_N(\theta)$.

It will be helpful to construct a POMP model $\bar f$ which expands the
model $f$ by allowing the parameter values to change deterministically
at each time point.
Specifically, we define a sequence $\vartheta_{0\dvtx N}=(\vartheta
_0,\ldots,\vartheta_N) \in\{\R^{\dimtheta}\}^{N+1}$. We then write
$(\barX_{0\dvtx N},\barY_{1\dvtx N})\in\RX^{N+1}\times\RY^N$ for a POMP
with generic density $\barf$ specified by the joint density
%
%
\begin{eqnarray}
&&\bar f_{\barX_{0\dvtx N},\barY_{1\dvtx N}}(x_{0\dvtx N},y_{1\dvtx N} ; \vartheta_{0\dvtx n})\nonumber
\\[-8pt]
\\[-8pt]
&& \qquad =
f_{X_0}(x_0;\vartheta_0)\prod_{k=1}^N f_{Y_k,X_k|
X_{k-1}}(y_k,x_k\given x_{k-1} ; \vartheta_k).
\nonumber
\end{eqnarray}
We write the log likelihood of $\data_{1\dvtx n}$ for the model $\barf$ as
$\barloglik_n(\vartheta_{0\dvtx n})=\break\log\bar f_{\barY_{1\dvtx n}}(\data
_{1\dvtx n};  \vartheta_{0\dvtx n})$.
We write $\theta^{[k]}$ to denote $k$ copies of $\theta\in\R
^\dimtheta$, concatenated in a column vector, so that $\loglik
_n(\theta)=\barloglik_n(\theta^{[n+1]})$.
We write $\nabla_i\barloglik_n(\vartheta_{0\dvtx n})$ for the partial
derivative of
$\barloglik_n$ with respect to $\vartheta_i$, for $i=0,\ldots,n$.
An application of the chain rule gives the identity
%
%
\begin{equation}\label{chainRule}
\nabla\loglik_n(\theta)
=\sum_{i=0}^n\nabla_i\barloglik_n\bigl(\theta^{[n+1]}\bigr).
\end{equation}
The regularity condition employed for Theorem~\ref{th1} below is
written in terms of this deterministically perturbed model:
\begin{longlist}
\item[(A4)]\hypertarget{A:loglik}
For each $1\le n\le N$, $\barloglik_n(\theta_{0\dvtx n})$ is twice differentiable.
For any compact subset $\compact$ of $\{\R^\dimtheta\}^{n+1}$ and
each $0\le i\le n$,
%
%
\begin{equation}
\sup_{\theta_{0\dvtx n}\in\compact}|\nabla_i\barloglik_n(\theta
_{0\dvtx n})|<\infty \quad
\mbox{and} \quad  \sup_{\theta_{0\dvtx n}\in\compact}|\nabla
_i^2\barloglik_n(\theta_{0\dvtx n})|<\infty.
\end{equation}
\end{longlist}
Condition \hyperlink{A:loglik}{(A4)} is a nonrestrictive
smoothness assumption.
However, the relationship between smoothness of the likelihood
function, the transition density $f_{X_k|X_{k-1}}(x_k\given
x_{k-1}; \theta)$, and the observation density $f_{Y_k|
X_k}(y_k\given x_k; \theta)$ is simple to establish only under the
restrictive condition that $\RX$ is a compact set.
Therefore, we note an alternative to \hyperlink{A:loglik}{(A4)} which
is more restrictive but more readily checkable:
\begin{longlist}
\item[(A4$'$)] $\RX$ is
compact. Both $f_{X_k|X_{k-1}}(x_k\given x_{k-1}; \theta)$
and $f_{Y_k|X_k}(y_k\given x_k; \theta)$ are twice
differentiable with respect to $\theta$.
These derivatives are continuous with respect to $x_{k-1}$ and $x_k$.
\end{longlist}

Iterated filtering involves introducing an auxiliary POMP model in
which a time-varying parameter process $\{\chTheta_n, 0 \le n
\le N\}$ is introduced.
Let $\kernel$ be a probability density function on $\R^\dimtheta$
having compact support, zero mean and covariance matrix $\ParScale$.
Let $Z_{0},\ldots,Z_N$ be $N$ independent draws from $\kernel$.
We introduce two perturbation parameters, $\sigma$ and $\tau$, and
construct a process $\chTheta_{0\dvtx N}$ by setting $\chTheta_0=\vartheta
_0+\tau Z_0$ and $\chTheta_k=\vartheta_k+\tau Z_0+\sigma\sum
_{j=1}^kZ_j$ for $1\leq k\leq N$.
The\vspace*{-1pt} joint density of $\chTheta_{0\dvtx N}$ is written as $g_{\chTheta
_{0\dvtx N}}(\breve{\vartheta}_{0\dvtx N}; \vartheta_{0\dvtx N},\sigma,\tau)$.
We define the\vspace*{-1pt} stochastically perturbed POMP model $g$ with a Markov
process $\{(\chX_n,\chTheta_n),0 {\le} n {\le} N\}$,
observation process $\chY_{1\dvtx N}$ and parameter $(\vartheta
_{0\dvtx N},\sigma,\tau)$ by the joint density
\begin{eqnarray*}
&&g_{\chX_{0\dvtx N},\chTheta_{0\dvtx N},\chY_{1\dvtx N}}(x_{0\dvtx N},\breve{\vartheta}
_{0\dvtx N},y_{1\dvtx N}; \vartheta_{0\dvtx N},\sigma,\tau)
\\
&& \qquad = g_{\chTheta_{0\dvtx N}}(\breve{\vartheta}_{0\dvtx N}; \vartheta
_{0\dvtx N},\sigma,\tau)
\bar f_{\barX_{0\dvtx N},\barY_{1\dvtx N}}(x_{0\dvtx N},y_{1\dvtx N} ; \breve
{\vartheta}_{0\dvtx N}).
\end{eqnarray*}
We seek a result analogous to Theorem~\ref{thm1} which takes into
account the specific structure of a POMP.
Theorem~\ref{th1} below gives a way to approximate $\nabla\loglik
_N(\theta)$ in terms of moments of the filtering distributions for $g$.
We write $\chE_{\vartheta_{0\dvtx n},\sigma,\tau}$ and $\chVar
_{\vartheta_{0\dvtx n},\sigma,\tau}$ for the expectation and variance,
respectively, for the model $g$.
We will be especially interested in the situation where $\vartheta
_{0\dvtx n}=\theta^{[n+1]}$, which leads us to define the following
filtering means and prediction variances:
%
%
\begin{eqnarray}\label{predmv}
\filttheta_n
&=&\filttheta_n(\theta,\sigma,\tau)=\chE_{\theta^{[n+1]},\sigma
,\tau}
[\chTheta_n\given\chY_{1\dvtx n}\equals\data_{1\dvtx n} ]\nonumber\\
 &=& \int
\breve{\vartheta}_n g^{}_{\chTheta_n|\chY_{1\dvtx n}}\bigl(\breve
{\vartheta}_n\given
\data_{1\dvtx n}; \theta^{[n+1]},\sigma,\tau\bigr) \, d\breve{\vartheta
}_n,\\
\predV_n &=& \predV_n(\theta,\sigma,\tau)=
\chVar_{\theta^{[n+1]},\sigma,\tau}
(\chTheta_n\given\chY_{1\dvtx n-1}\equals\data_{1\dvtx n-1} )
\nonumber
\end{eqnarray}
for $n=1,\ldots,N$, with $\filttheta_0=\theta$.

%
\begin{theorem}\label{th1}
Suppose condition \hyperlink{A:loglik}{(\textup{A}4)}.
Let $\sigma$ be a function of $\tau$ with $\lim_{\tau\to0}\sigma
(\tau) /\tau=0$.
For any compact set $\kD\subset\R^\dimtheta$, there exists a finite
constant $\cLL$ such that for all $\tau$ small enough,
%
%
\begin{equation}
\sup_{\theta\in\kD} |
\tau^{-2} \Sigma^{-1} (\filttheta_N-\theta)
- \nabla\loglik_N(\theta) |
\le\cLL(\tau+\sigma^2/\tau^2)
\label{eqT101}
\end{equation}
and
%
%
\begin{equation}
\sup_{\theta\in\kD} \Biggl|
\sum_{n=1}^N ( \predV_n )^{-1}
(
\filttheta_n-\filttheta_{n-1}
)
-\nabla\loglik_N(\theta)
\Biggr|
\le\cLL(\tau+\sigma^2/\tau^2).
\label{eqT1}
\end{equation}
\end{theorem}

\begin{pf}
For each $n\in\{1,\ldots,N\}$, we map onto the notation of
Section~\ref{seciis} by setting $X=\barX_{0\dvtx n}$, $Y=\barY_{1\dvtx n}$,
$\theta=\vartheta_{0\dvtx n}$, $\chTheta=\chTheta_{0\dvtx n}$, $\data=\data
_{1\dvtx n}$ and $h(\vartheta_{0\dvtx n})=\vartheta_n$.
We note that, by construction, this implies $\chX=\chX_{0\dvtx n}$, $\chY
=\chY_{1\dvtx n}$, and $\kernel_\tau(\vartheta_{0\dvtx n})=\kernel(\vartheta
_0)\prod_{i=1}^n(\sigma/\tau)^{-\dimtheta}\kernel((\vartheta
_i-\vartheta_{i-1})/(\sigma/\tau) )$.
For this choice of $h$, the integral $\tau\int h(\tau u)u'\kernel
_\tau(u)\, du$ is a $p\times p(n+1)$ matrix for which the $i$th
$p\times p$ sub-matrix is
$\chE_{\vartheta_{0\dvtx n},\sigma,\tau} [(\chTheta_n-\vartheta
_n)(\chTheta_{i-1}-\vartheta_{i-1})^\transpose]=\{\tau
^2+(i-1)\sigma^2\}\Sigma$.
Thus,
%
%
\begin{eqnarray}
\label{eqterm3}
&&\biggl(\tau\int h(\tau u)u'\kernel_\tau(u)\, du \biggr)\nabla
\bar\ell_n(\vartheta_{0\dvtx n})\nonumber\\
&& \qquad =
\sum_{i=0}^n\chE_{\vartheta_{0\dvtx n},\sigma,\tau}
[
(\chTheta_n-\vartheta_n)(\chTheta_i-\vartheta_i)^\transpose
]
\nabla_i\bar\ell_n(\vartheta_{0\dvtx n})
\\
&& \qquad =
\sum_{i=0}^n (\tau^2\Sigma+i \sigma^2 \Sigma)\nabla
_i\bar\ell_n(\vartheta_{0\dvtx n}).
\nonumber
\end{eqnarray}
Applying Theorem~\ref{thm1} in this context, the second term in (\ref{eq555}) is zero and the third term is given by
(\ref{eqterm3}).
We obtain that for\ any compact $\kH\subset\R^{p(n+1)}$ there is a
$\cMM<\infty$ such that
%
%
\begin{eqnarray}
\label{eq1thmIF}
&&\sup_{\vartheta_{0\dvtx n}\in\kH}
\Biggl|\tau^{-2}\Sigma^{-1} (\chE_{\vartheta_{0\dvtx n},\sigma
,\tau}[\chTheta_n \given Y_{1\dvtx n}=y_{1\dvtx n}^*]-\vartheta_n
)\nonumber\\
&&\hspace*{90pt}{}- \sum_{i=0}^n (1+\sigma^2\tau^{-2} i )\nabla_i\bar\ell
_n(\vartheta_{0\dvtx n})
\Biggr|\\
&& \qquad
\leq\cMM\tau.\nonumber
\end{eqnarray}
Applying (\ref{eq1thmIF}) to the special case of $\vartheta
_{0\dvtx n}=\theta^{[n+1]}$, making use of (\ref{chainRule}) and (\ref{predmv}), we infer the existence of finite constants $\cMMM$ and
$\cNN$ such that
%
%
\begin{eqnarray}
\label{eq1thmIFv2}
 \qquad \sup_{\theta\in\kD} |\tau^{-2}\Sigma^{-1} (\filttheta
_n-\theta) -\nabla\ell_n(\theta) |
&\leq&
\cMMM\tau
+\sup_{\theta\in\kD}\sigma^2\tau^{-2}\sum_{i=1}^ni \bigl|\nabla
_i\ell_n\bigl(\theta^{[n+1]}\bigr) \bigr|
\nonumber
\\[-8pt]
\\[-8pt]
 \qquad &\leq&
\cNN(\tau+\sigma^2\tau^{-2}),
\nonumber
\end{eqnarray}
which establishes (\ref{eqT101}).
To show (\ref{eqT1}), we write
%
%
\begin{eqnarray}
\label{varIdentity}
&&\sum_{k=1}^n \{\predV_k \}^{-1} (\filttheta_k-\filttheta
_{k-1} )\nonumber\\
&& \qquad =
\sum_{k=1}^n \tau^{-2}\Sigma^{-1} (\filttheta_k-\filttheta
_{k-1} )
+\tau^{-2}\sum_{k=1}^n
(
\{\tau^{-2}\predV_k \}^{-1}
-
\Sigma^{-1}
)
(
\filttheta_k-\filttheta_{k-1}
)
\\
&& \qquad =
\tau^{-2}\Sigma^{-1}
(
\filttheta_n-\theta
)
+ \tau^{-2}\sum_{k=1}^n
( \{\tau^{-2}\predV_k \}^{-1} -
\Sigma^{-1}
)
(\filttheta_k-\filttheta_{k-1} ).
\nonumber
\end{eqnarray}
We note that (\ref{eq1thmIFv2}) implies the existence of a bound
%
%
\begin{equation} \label{filtThetaBound}
|\filttheta_n-\theta|\le\cQQ\tau^2.
\end{equation}
Combining (\ref{varIdentity}), (\ref{eq1thmIFv2}) and (\ref{filtThetaBound}), we deduce
\begin{eqnarray*}
&&\sup_{\theta\in\kD} \Biggl|\sum_{k=1}^n \{\predV_k \}
^{-1} (\filttheta_{k}-\filttheta_{k-1} ) -\nabla\ell
_n(\theta) \Biggr|\\
&& \qquad \leq\cOO(\tau+\sigma^2\tau^{-2} ) +
\cRR\sum_{k=1}^n\sup_{\theta\in\kD}
|
\{\tau^{-2}\predV_k \}^{-1} -
\Sigma^{-1}
|
\end{eqnarray*}
for finite constants $\cOO$ and $\cRR$.
For invertible matrices $A$ and $B$, we have the bound
%
%
\begin{equation}\label{matid}
|A^{-1}-B^{-1} |\leq|B^{-1} |^2
\bigl(
1- |(B-A)B^{-1} |
\bigr)^{-1}
|B-A |
\end{equation}
provided that $ |(B-A)B^{-1} |<1$.
Applying (\ref{matid}) with $A= \{\tau^{-2}\predV_k \}$ and
$B=\Sigma$,
we see that the theorem will be proved once it is shown that
\[
\sup_{\theta\in\kD} |
\tau^{-2}\predV_k -
\Sigma|
\leq\cPP(\tau+\sigma^2\tau^{-2} ).
\]
Now, it is easy to check that
\begin{eqnarray*}
\tau^{-2}\predV_n-\Sigma&=&\tau^{-2}\chE_{\theta^{[n+1]},\sigma
,\tau}
[
(\chTheta_{n-1}-\theta) (\chTheta_{n-1}-\theta
)^\transpose
\mid \chY_{1\dvtx n-1}\equals\data_{1\dvtx n-1}
]-\Sigma\\
&&{}
-\tau^{-2} (\filttheta_{n-1}-\theta) (\filttheta
_{n-1}-\theta)^\transpose+ \sigma^2\tau^{-2}\Sigma.
\end{eqnarray*}
Applying Theorem~\ref{thm1} again with $h(\vartheta_{0\dvtx n})=(\vartheta
_n-\theta)(\vartheta_n-\theta)^\transpose$, and making use of (\ref{filtThetaBound}), we obtain
%
%
\begin{eqnarray}\label{eq2thmIF}
&&\sup_{\theta\in\kD} \bigl|\tau^{-2}\chE_{\theta^{[n+1]},\sigma
,\tau}
[
(\chTheta_n-\theta) (\chTheta_n-\theta
)^\transpose
\mid
\chY_{1\dvtx n}\equals\data_{1\dvtx n}
]
-\Sigma\bigr|\nonumber
\\[-8pt]
\\[-8pt]
&& \qquad \leq\cSS(\tau+\sigma^2\tau^{-2} ),
\nonumber
\end{eqnarray}
which completes the proof.
\end{pf}

The two approximations to the derivative of the log likelihood in (\ref{eqT101}) and (\ref{eqT1})
are asymptotically equivalent in the
theoretical framework of this paper.
However, numerical considerations may explain why (\ref{eqT1}) has
been preferred in practical applications.
To be concrete, we suppose henceforth that numerical filtering will be
carried out using the basic sequential Monte Carlo method presented as
Algorithm~\ref{algpf}.
Sequential Monte Carlo provides a flexible and widely used class of
filtering algorithms, with many variants designed to improve numerical
efficiency [\citet{cappe07}].
The relatively simple sequential Monte Carlo method in Algorithm~\ref{algpf} has,
however, been found adequate for previous data analyses
using iterated filtering [\citet{ionides06-pnas};
\citet{king08};
\citet{breto09};
\citet{he10};
\citet{laneri10}].

%
%
\begin{algorithm}[t]
\def\qsp{\hspace{0.5cm}}
\textbf{Input:}
\begin{itemize}
\item POMP model described by a generic density $\genf$ having
parameter vector $\gentheta$ and corresponding to a Markov process
$\genX_{0\dvtx N}$, observation process $\genY_{1\dvtx N}$, and data $\gendata_{1\dvtx N}$
\item
Number of particles, $J$
\end{itemize}

\textbf{Procedure:} 
\begin{enumerate}[3]
\item[1] initialize filter particles $\genX^\filter_{0,j}\sim\genf
_{\genX_0}(\tilde{x}_0 ; \gentheta)$ for
$j$ in $1\dvtx J$
\item[2] for $n$ in $1\dvtx N$
\item[3]\hypertarget{step:pred}
\qsp for $j$ in $1\dvtx J$ 
draw prediction particles $\genX^\predict_{n,j}\sim\genf_{\genX
_n|\genX_{n-1}} (\tilde{x}_n\given\genX^\filter_{{n-1},j}
; \gentheta)$
\item[4]\hypertarget{step:w}
\qsp set $w(n,j)=\genf_{\genY_n|\genX_n} (\gendata
_n\given\genX^\predict_{{n},j} ; \gentheta)$\vspace*{1pt}
\item[5] \hypertarget{step:sample}\qsp draw $k_1,\ldots,k_{J}$ such
that $\prob\{k_j{=}i\}
=w(n,i)/\sum_\ell w(n,\ell)$
\item[6]\hypertarget{step:filter}
\qsp set $\genX^\filter_{{n},j}=\genX^\predict_{{n},k_j}$
\item[7] end for
\end{enumerate}
\caption{A basic sequential Monte Carlo procedure for a discrete-time
Markov process.
For the unperturbed model, set $\genX_n=X_n$, $\genY_n=Y_n$,$\genf
=f$ and $\gentheta=\theta$.
For the stochastically perturbed model, set $\genX_n=(\chX_n,\chTheta
_n)$, $\genY_n=Y_n$, $\genf=g$ and $\gentheta=(\theta
^{[N+1]},\sigma,\tau)$.
It is neither necessary nor computationally optimal to draw from the
density $\genf_{\genX_n|\genX_{n-1}}$ in step~\protect\hyperlink
{step:pred}{3} [e.g., \protect\citet{arulampalam02}] however only this choice
leads to the plug-and-play property.
The resampling in step~\protect\hyperlink{step:sample}{5} is taken to follow a
multinomial distribution to build on previous theoretical results
making this assumption [\protect\citet{delmoral01};
\protect\citet{crisan02}].
An alternative is the systematic procedure in \protect\citeauthor{arulampalam02}
[(\protect\citeyear{arulampalam02}), Algorithm~2] which has less Monte Carlo
variability.
We support the use of systematic sampling in practice, and we suppose
that all our results would continue to hold in such situations.}
\label{algpf}
\end{algorithm}

When carrying out filtering via sequential Monte Carlo, the resampling
involved has a consequence that all surviving particles can descend
from only few recent ancestors.
This phenomenon, together with the resulting shortage of diversity in
the Monte Carlo sample, is called \textit{particle depletion} and can
be a
major obstacle for the implementation of sequential Monte Carlo
techniques [\citet{arulampalam02};
\citet{andrieu10}].
The role of the added variation on the scale of $\sigma_m$ in the
iterated filtering algorithm is to rediversify the particles and hence
to combat particle depletion.
Mixing considerations suggest that the new information about $\theta$
in the $n$th observation may depend only weakly on $\data_{1\dvtx n-k}$ for
sufficiently large~$k$ [\citet{jensen99b}].
The actual Monte Carlo particle diversity of the filtering
distribution, based on (\ref{eqT1}), may therefore be the best guide
when sequentially estimating the derivative of the log likelihood.
Future theoretical work on iterated filtering algorithms should study
formally the role of mixing, to investigate this heuristic argument.
However, the theory presented in Theorems~\ref{th1},~\ref{th2}
and~\ref{th3} does formally support using a limit random walk
perturbations without any mixing conditions.
Two influential previous proposals to use stochastic perturbations to
reduce numerical instabilities arising in plug-and-play inference for
POMPS via sequential Monte Carlo [\citet{kitagawa98};
\citet{janeliu01}] lack
even this level of theoretical support.

To calculate Monte Carlo estimates of the quantities in (\ref{predmv}),
 we\vspace*{1pt} apply Algorithm~\ref{algpf} with $\genf=g$, $\genX
_n=(\chX_n,\chTheta_n)$, $\gentheta=(\theta^{[N+1]},\sigma,\tau)$
and $J$ particles.\vspace*{1pt}
We write $\genX^\predict_{n,j}= ( X^\predict_{n,j} ,\Theta
^\predict_{n,j} )$ and $\genX^\filter_{n,j}= ( X^\filter
_{n,j} ,\Theta^\filter_{n,j} )$ for the Monte Carlo\vspace*{1pt} samples
from the prediction and filtering and calculations in steps~\hyperlink
{step:pred}{3} and~\hyperlink{step:filter}{6} of Algorithm~\ref{algpf}.
Then we define
%
%
\begin{eqnarray}\label{predmvmc}
\widetilde{\theta}^\filter_n&=&\widetilde\theta^\filter_n(\theta
,\sigma,\tau,J)=
\frac{1}{J} \sum_{j=1}^J \Theta^\filter_{n,j}, \nonumber
\\[-8pt]
\\[-8pt]
\widetilde V^\predict_n&=&\widetilde V^\predict_n(\theta,\sigma,\tau,J)=
\frac{1}{J-1} \sum_{j=1}^J
( \Theta^\predict_{n,j} - \widetilde\theta^\filter_{n-1} )
( \Theta^\predict_{n,j} - \widetilde\theta^\filter_{n-1}
)^\transpose
.
\nonumber
\end{eqnarray}
In practice, a reduction in Monte Carlo variability is possible by
modifying (\ref{predmvmc}) to estimate ${\theta}^\filter_n$ and
$V^\predict_n$ from weighted particles prior to resampling
[\citet{chopin04}].
We now present, as Theorem~\ref{th2}, an analogue to Theorem~\ref{th1} in which the filtering means and prediction variances are
replaced by their Monte Carlo counterparts.
The stochasticity in Theorem~\ref{th2} is due to Monte Carlo
variability, conditional on the data $\data_{1\dvtx N}$, and we write $\EMC
$ and $\VarMC$ to denote Monte Carlo means and variances.
The Monte Carlo random variables required to implement Algorithm~\ref{algpf} are presumed to be drawn independently each time the algorithm
is evaluated.
To control the Monte Carlo bias and variance, we assume:
\begin{longlist}
\item[(A5)]\hypertarget{A:fbound}
For each $n$ and any compact set $\kE\subset\R^\dimtheta$,
\[
\sup_{\theta\in\kE, x\in\RX} f_{Y_n|X_n}(\data
_n\given x_n; \theta)<\infty.
\]
\end{longlist}
%
%
\begin{theorem}\label{th2} 
Let $\{\sigma_m\}$, $\{\tau_m\}$ and $\{J_m\}$ be positive sequences
with \mbox{$\tau_m\to0$}, $\sigma_m \tau^{-1}_m\to0$ and $\tau_m J_m\to
\infty$.
Define $\widetilde\theta^\filter_{n,m}=\widetilde\theta^\filter
_n(\theta,\sigma_m,\tau_m,J_m)$ and
$\widetilde V^\predict_{n,m}=\widetilde V^\predict_{n}(\theta,\sigma
_m,\tau_m,J_m)$ via (\ref{predmvmc}).
Suppose conditions \hyperlink{A:loglik}{(\textup{A}4)} and \hyperlink
{A:fbound}{(\textup{A}5)} and let $\kF$ be an arbitrary compact subset
of $\R^\dimtheta$. Then,
%
%
\begin{eqnarray} \label{eqT2a}
\lim_{m\rightarrow\infty}\sup_{\theta\in\kF} \Biggl|\EMC\Biggl[
\sum_{n=1}^N (\widetilde V^\predict_{n,m} )^{-1}
(
\widetilde\theta^\filter_{n,m}-\widetilde\theta^\filter_{n-1,m}
)
\Biggr]
-
\nabla\loglik_N(\theta)
\Biggr|
&=& 0
,\\ \label{eqT2b}
\lim_{m\rightarrow\infty}
\sup_{\theta\in\kF}
\Biggl|
\tau_m^2 J_m
\VarMC\Biggl( \sum_{n=1}^N (\widetilde V^\predict_{n,m} )^{-1}
(
\widetilde\theta^\filter_{n,m}-\widetilde\theta^\filter_{n-1,m}
) \Biggr)
\Biggr|
&<&
\infty.
\end{eqnarray}
%
\end{theorem}

\begin{pf}
Let $\kI$ be a compact subset of $\R^2$ containing $\{(\sigma_m,\tau
_m), m=1,2,\ldots\}$.
Set $\theta\in\kF$ and $(\sigma,\tau)\in\kI$.
Making use of the definitions in (\ref{predmv}) and (\ref{predmvmc}), we construct $u_{n}=(\filttheta_{n}-\filttheta
_{n-1})/\tau$ and $v_{n}= \predV_{n}/\tau^2$, with corresponding
Monte Carlo estimates $\widetilde u_{n}=(\widetilde\theta^\filter
_{n}-\widetilde\theta^\filter_{n-1})/\tau$ and $\widetilde
v_{n}=\widetilde V^\predict_{n}/\tau^2$.\vspace*{-1pt}
We look to apply Theorem~\ref{thcrisan} (presented in Section~\ref{appsmc}) with\vspace*{1pt} $\genf=g$, $\genX_n=(\chX_n,\chTheta_n)$, $\genY
_n=\chY_n$, $\gentheta=(\theta^{[n+1]},\sigma,\tau)$, $J$
particles, and
\[
\phi(\chX_n,\chTheta_n)= (\chTheta_n -\theta)  / \tau.
\]
Using the notation from (\ref{eqphi}), we have
$u_{n}=\testFunction_n^F-\testFunction_{n-1}^F$
and
$\widetilde u_{n}=\widetilde\testFunction^F_n-\widetilde
\testFunction^F_{n-1}$.
By assumption, $\kernel(u)$ is supported on some set $\{u\dvtx  |u|<\cD\}$
from which we derive the bound $ | \phi(\check X_n, \check\Theta
_n) | \le\cD(1+n\sigma/\tau)$.
Theorem~\ref{thcrisan} then provides for the existence of a $\cW$
and $\cY$ such that
%
%
\begin{eqnarray}
\label{eq501}
\EMC[ |\widetilde u_{n} - u_{n}|^2 ] &\le& \cW/J,
\\
\label{eq502}
| \EMC[ \widetilde u_{n} - u_{n} ] | &\le&
\cY/J.
\end{eqnarray}
The explicit bounds in (\ref{eqgamma}) of Theorem~\ref{thcrisan},
together with \hyperlink{A:loglik}{(A4)} and \hyperlink
{A:fbound}{(A5)}, assure us that $\cW=\cW(\theta,\sigma,\tau)$
and $\cY=\cY(\theta,\sigma,\tau)$ can be chosen so that (\ref{eq501}) and (\ref{eq502}) hold uniformly over $(\theta,\sigma
,\tau)\in\kF\times\kI$.
The same argument applied to $v_{n}= V^\predict_{n}/\tau^2$ and
$\widetilde v_{n}=\widetilde V^\predict_{n}/\tau^2$ gives
%
%
\begin{equation}
| \EMC[\widetilde v_{n} - v_{n} ] | \le\cZ
/J , \qquad
\EMC[ |\widetilde v_{n} - v_{n}|^2 ] \le\cB/J
\label{eq503}
\end{equation}
uniformly over $(\theta,\sigma,\tau)\in\kF\times\kI$.
We now proceed to carry out a Taylor series expansion:
%
%
\begin{equation}
{\widetilde v_{n}^{-1}}\widetilde u_{n}
=
{v_{n}^{-1}}
u_{n}
+
{v_{n}^{-1}} (\widetilde u_{n} -u_{n})
{v_{n}^{-1}}(\widetilde v_{n} -v_{n}){v_{n}^{-1}} u_{n}
+\rB,\label{eqtaylor1}
\end{equation}
where $|\rB| < \cV(|\widetilde u_{n} -u_{n}|^2 + |\widetilde v_{n}
-v_{n}|^2)$ for some constant $\cV$.
The existence of such a $\cV$ is guaranteed since the determinant of
$v_{n}$ is bounded away from zero.
Taking expectations of both sides of (\ref{eqtaylor1}) and applying
(\ref{eq501})--(\ref{eq503}) gives
%
%
\begin{equation}
|
\EMC[{\widetilde v_{n}^{-1}}\widetilde u_{n}] -
v_{n}^{-1} u_{n}
|
\le\cT/J
\label{eqtaylor2}
\end{equation}
for some constant $\cT<\infty$. Another Taylor series expansion,
\[
{\widetilde v_{n}^{-1}}\widetilde u_{n}
=
{v_{n}^{-1}}
u_{n}
+\rC
\]
with $|\rC| < \cU(|\widetilde u_{n} -u_{n}| + |\widetilde v_{n}
-v_{n}|)$ implies
%
%
\begin{equation}
\VarMC(\widetilde v_{n}^{-1}\widetilde u_{n})
\le\cS/J.
\label{eqtaylor4}
\end{equation}
Rewriting (\ref{eqtaylor2}) and (\ref{eqtaylor4}), defining
$\filttheta_{n,m}=\filttheta_n(\theta,\sigma_m,\tau_m)$ and
$\predV_{n,m}=\predV_n(\theta,\break\sigma_m, \tau_m)$,
we deduce that
%
%
\begin{equation}
\label{eqfiltMeanDiff}
 \qquad \tau_mJ_m |
\EMC
[ (\widetilde V^\predict_{n,m} )^{-1}(\widetilde\theta
^\filter_{n,m}-\widetilde\theta^\filter_{n-1,m})
]
- (\predV_{n,m} )^{-1}(\filttheta_{n,m}-\filttheta_{n-1,m})
|
\le\cT
\end{equation}
and
%
%
\begin{equation}
\label{eqfiltVarDiff}
\tau_m^2 J_m
\VarMC[
(\widetilde V^\predict_{n,m} )^{-1}(\widetilde\theta
^\filter_{n,m}-\widetilde\theta^\filter_{n-1,m})
] \le\cS.
\end{equation}
Combining (\ref{eqfiltMeanDiff}) with Theorem~\ref{th1}, and
summing over $n$, leads to (\ref{eqT2a}).
Summing (\ref{eqfiltVarDiff}) over $n$ justifies (\ref{eqT2a}).
\end{pf}

Theorem~\ref{th2} suggests that a Monte Carlo method which leans on
Theorem~\ref{th1} will require a sequence of Monte Carlo sample
sizes, $J_m$, which increases faster than $\tau_m^{-1}$.
Even with $\tau_m J_m\to\infty$, we see from (\ref{eqT2b}) that
the estimated derivative in (\ref{eqT2a}) may have increasing
Monte Carlo variability as $m\to\infty$.
Theorem~\ref{th3} gives an example of a stochastic approximation
procedure, defined by the recursive sequence $\hat\theta_m$ in (\ref{eqrecursion}), that makes use of the Monte Carlo estimates studied in
Theorem~\ref{th2}.
Because each step of this recursion involves an application of the
filtering procedure in Algorithm~\ref{algpf}, we call (\ref{eqrecursion}) an iterated filtering algorithm.
The rate assumptions in Theorem~\ref{th3} are satisfied, for example,
by $a_m = m^{-1}$, $\tau_m^2 = m^{-1}$, $\sigma^2_m=m^{-(1+\delta)}$
and $J_m = m^{(\delta+1/2)}$ for $\delta>0$.

%
%
\begin{theorem}\label{th3} 
Let $\{a_m\}$, $\{\sigma_m\}$, $\{\tau_m\}$ and $\{J_m\}$ be
positive sequences with \mbox{$\tau_m\to0$}, $\sigma_m\tau_m^{-1}\to0$,
$ J_m\tau_m \to\infty$, $a_m\to0$, $\sum_m a_m=\infty$ and $\sum
_m a_m^2J_m^{-1}\tau_m^{-2} <\infty$.
Specify a recursive sequence of parameter estimates $\{\hat\theta_m\}
$ by
%
%
\begin{equation}
\hat\theta_{m+1} = \hat\theta_m + a_m \sum_{n=1}^N
(\widetilde V^\predict_{n,m} )^{-1}
(
\widetilde\theta^\filter_{n,m}-\widetilde\theta^\filter_{n-1,m}
), \label{eqrecursion}
\end{equation}
where $\widetilde\theta^\filter_{n,m}=\widetilde\theta^\filter
_n(\hat\theta_{m},\sigma_m,\tau_m,J_m)$ and
$\widetilde V^\predict_{n,m}=\widetilde V^\predict_{n,m}(\hat\theta
_{m},\sigma_m,\tau_m,J_m)$ are defined in (\ref{predmvmc}) via an
application of Algorithm~\ref{algpf}.
Assuming \hyperlink{A:loglik}{(\textup{A}4)}, \hyperlink{C1}{(\textup{B}1)} and
\hyperlink{C2}{(\textup{B}2)}, $\lim_{m\to\infty} \hat
\theta_m=\hat\theta$ with probability one.
\end{theorem}

\begin{pf}
Theorem~\ref{th3} follows directly from a general stochastic
approximation result, Theorem~\ref{thspall} of Section~\ref{appsa}, applied to $\loglik_N(\theta)$.
Conditions \hyperlink{C4}{(B4)} and \hyperlink{C5}{(B5)} of
Theorem~\ref{thspall} hold from Theorem~\ref{th2} and the
remaining assumptions of Theorem~\ref{thspall} hold by hypothesis.
\end{pf}

\section{Discussion}
\label{secdiscussion} 

One alternative approach to likelihood maximization for POMP models
involves plugging the (log) likelihood estimate from a particle filter
directly into a general-purpose stochastic optimization algorithm such
as Simultaneous Perturbation Stochastic Approximation (SPSA),
Kiefer--Wolfowitz or stochastic Nelder--Mead [\citet{spall03}].
An advantage of iterated filtering, and other methods based on particle
filtering with parameter perturbations [\citet{kitagawa98};
\citet{janeliu01}],
is that the many thousands of particles are simultaneously exploring
the parameter space and evaluating an approximation to the likelihood.
When the data are a long time series, the perturbed parameters can make
substantial progress toward plausible parameter values in the course of
one filtering operation.
From the point of view of a general-purpose stochastic optimization
algorithm, carrying out one filtering operation (which can be a
significant computational burden in many practical situations) yields
only one function evaluation of the likelihood.

The practical applicability of particle filters may be explained by
their numerical stability on models possessing a mixing property [e.g.,
\citet{crisan02}].
The sequential Monte Carlo analysis in Theorem~\ref{th2} did not
address the convergence of iterated filtering under mixing assumptions
as the number of observations, $N$, increases.
We therefore studied experimentally the numerical stability of the
Monte Carlo estimate of the derivative of the log likelihood in
equation (\ref{eqT2a}).
The role of mixing arises regardless of the dimension of the state
space, the dimension of the parameter space, the nonlinearity of the
system, or the non-Gaussianity of the system.
This suggests that a simple linear Gaussian example may be
representative of behavior on more complex models.
Specifically, we considered a POMP model defined by a scalar Markov
process $X_n=\theta X_{n-1}+\epsilon_n$, with $X_0=0$, and a scalar
observation process $Y_n=X_n+\eta_n$.
Here, $\{\epsilon_n\}$ and $\{\eta_n\}$ were taken to be sequences of
independent Gaussian random variables having zero mean and unit variance.
We fixed the true parameter value as $\theta^*=0.8$ and we evaluated
$\nabla\loglik_N(\theta)$ at $\theta=\theta^*$ and $\theta=0.9$
using a Kalman filter (followed by a finite difference derivative
computation) and via the sequential Monte Carlo approximation in (\ref{eqT2a}) using $J=1\mbox{,}000$ particles.
We investigated $\sigma\in\{0.002,0.005,0.02\}$, chosen to include a
small value where Monte Carlo variance dominates, a large value where
bias dominates, and an intermediate value; we then fixed $\tau=20
\sigma$.

%
%
\begin{figure}

\includegraphics{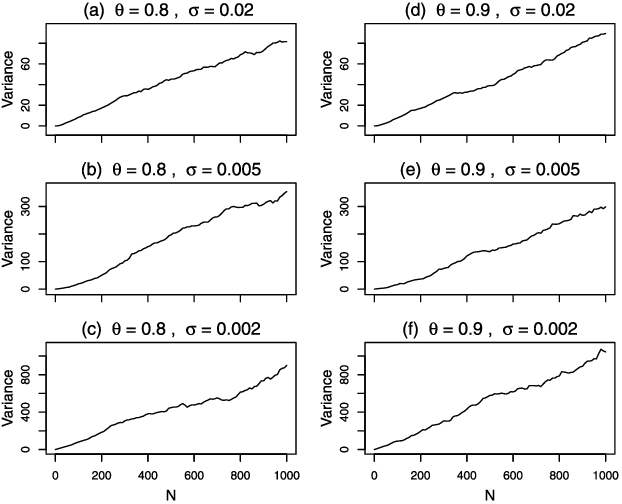}

\caption{Monte Carlo variance of the derivative approximation in
(\protect\ref{eqT2a}) for varying values of $\theta$ and~$\sigma$.}\vspace*{-6pt}\label{figvar}
\end{figure}

%
%
\begin{figure}

\includegraphics{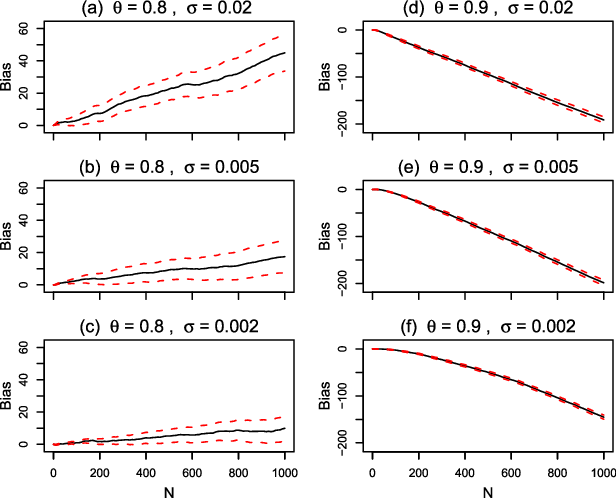}

\caption{Bias of the derivative approximation in (\protect\ref{eqT2a}) for
varying values of $\theta$ and $\sigma$. Dashed lines show pointwise
95\% Monte Carlo confidence intervals.}\label{figbias}\vspace*{-6pt}
\end{figure}

%
%
\begin{figure}
\includegraphics{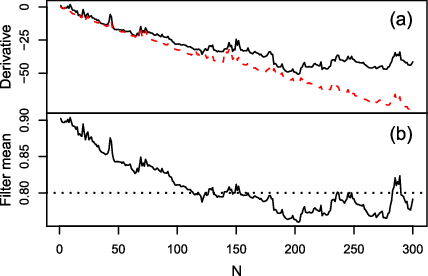}

\caption{One realization from the simulation study, with $\theta=0.9$
and $\sigma=0.005$. \textup{(a)} The estimate of $\nabla\loglik_N(\theta)$
using (\protect\ref{eqT2a}) (solid line) and calculated directly (dashed
line). \textup{(b)} The filter mean $\tilde\theta^F_N$ (solid line)
approaching the vicinity of the true parameter value $\theta^*=0.8$
(dotted line).}\label{figsim}
\end{figure}

Figures~\ref{figvar} and~\ref{figbias} show how the Monte Carlo
variance and the bias vary with $N$ for each value of $\sigma$.
These quantities were evaluated from 100 realizations of the model,
with 5 replications of the filtering operation per realization, via
standard unbiased estimators.
We see from Figure~\ref{figvar} that the Monte Carlo variance
increases approximately linearly with $N$.
This numerical stability is a substantial improvement on the
exponential bound guaranteed by Theorem~\ref{thcrisan}.
The ordinate values in Figure~\ref{figvar} show that, as anticipated
from Theorem~\ref{th2}, the variance increases as $\sigma$ decreases.
Figure~\ref{figbias} shows that the bias diminishes as $\sigma$
decreases and is small when $\theta$ is close to $\theta^*$.
When $\theta$ is distant from $\theta^*$, the perturbed parameter
values migrate toward $\theta^*$ during the course of the filtering
operation, as shown in Figure~\ref{figsim}(b).
Once the perturbed parameters have arrived in the vicinity of $\theta
^*$, the sum in (\ref{eqT2a}) approximates the derivative of the log
likelihood at $\theta^*$ rather than at $\theta$.
Figure~\ref{figsim}(a) demonstrates the resulting bias in the
estimate of $\nabla\loglik_N(\theta)$.
However, this bias may be helpful, rather than problematic, for the
convergence of the iterated filtering algorithm.
The update in (\ref{eqrecursion}) is a weighted average of the
filtered means of the perturbed parameters.
Heuristically, if the perturbed parameters successfully locate a
neighborhood of $\theta^*$ then this will help to generate a good
update for the iterated filtering algorithm.
The utility of perturbed parameter values to identify a neighborhood of
$\theta^*$, in addition to estimating a derivative, does not play a
role in our asymptotic justification of iterated filtering.
However, it may contribute to the nonasymptotic properties of the
method at early iterations.
%



\begin{appendix} 
\section*{Appendix: Some standard results on sequential Monte Carlo and
stochastic approximation theory}

We state some basic theorems that we use to prove Theorems~\ref{thsaiis},~\ref{th2} and~\ref{th3}, both for completeness and
because we require minor modifications of the standard results.\vadjust{\goodbreak}
Our goal is not to employ the most recent results available in these
research areas, but rather to show that some fundamental and well-known
results from both areas can be combined with our Theorems~\ref{thm1}
and~\ref{th1} to synthesize a new theoretical understanding of
iterated filtering and iterated importance sampling.


\subsection{A version of a standard stochastic approximation theorem}
\label{appsa}

We present, as Theorem~\ref{thspall}, a special case of Theorem~2.3.1
of \citet{kushner78}.
For variations and developments on this result, we refer the reader to
\citet{kushner03}, \citet{spall03}, \citet{andrieu05}
and \citet{maryak08}.
In particular, Theorem~2.3.1 of \citet{kushner78} is similar to
Theorem~4.1 of \citet{spall03} and to Theorem~2.1 of \citet
{kushner03}.

\begin{theorem}\label{thspall}
%
%
Let $\loglik(\theta)$ be a continuously differentiable function $\R
^\dimtheta\to\R$ and let $\{D_m(\theta), m\ge1\}$ be a sequence of
independent Monte Carlo estimators of the vector of partial derivatives
$\nabla\loglik(\theta)$.
Define a sequence $\{\hat\theta_m\}$ recursively by
$
\hat\theta_{m+1}=\hat\theta_{m} + a_mD_m(\hat\theta_{m})
$.
Assume \hyperlink{C1}{(\textup{B}1)} and \hyperlink{C2}{(\textup{B}2)} of
Section~\ref{seciis} together with the following conditions:
\begin{longlist}[(B5)]
\item[(B3)]
\hypertarget{C3}
$a_m>0$, $a_m\to0$, $\sum_m a_m=\infty$.
\item[(B4)]\hypertarget{C4}
$\sum_{m}a_m^2\sup_{|\theta|<r}\VarMC(D_m(\theta) )
<\infty$ for every $r>0$.
\item[(B5)]
\hypertarget{C5}
$\lim_{m\to\infty}\sup_{|\theta|<r}
|
\EMC[D_{m}(\theta)] - \nabla\loglik(\theta)
|
= 0$
for every $r>0$.
\end{longlist}
Then $\hat\theta_{m}$ converges to $\hat\theta=\arg\max\loglik
(\theta)$ with probability one.
\end{theorem}

\begin{pf}
The most laborious step in deducing Theorem~\ref{thspall} from
\citet
{kushner78} is to check that \hyperlink{C1}{(B1)}--\hyperlink
{C5}{(B5)} imply that, for all $\epsilon>0$,
%
%
\begin{equation}
\lim_{n\to\infty} \prob\Biggl[\sup_{j\ge1} \Biggl|
\sum_{m=n}^{n+j} a_m
\{
D_m(\hat\theta_{m})-\EMC[D_m(\hat\theta_{m})\given\hat\theta_{m}]
\}
\Biggr|
\ge\epsilon
\Biggr]=0,
\label{eqkc}
\end{equation}
which in turn implies condition A2.2.4$^{\prime\prime}$ and hence
A2.2.4 of \citet{kushner78}.
To show (\ref{eqkc}), we define
$\xi_m=D_m(\hat\theta_{m})-\EMC[D_m(\hat\theta_{m})\given\hat
\theta_{m}]$ and
%
%
\begin{equation}
\xi_m^k=
\cases{\displaystyle
\xi_m,
& \quad if $|\hat\theta_{m}|\le k$,\vspace*{1pt}\cr\displaystyle
0, & \quad  if $|\hat\theta_{m}|>k$.
}
\end{equation}
Define processes
$\{M^{n}_j{=}\sum_{m=n}^{n+j} a_m\xi_m, j\ge0\}$ and
$\{M^{n,k}_j{=}\sum_{m=n}^{n+j} a_m\xi_m^k, j\ge0\}$ for each $k$
and $n$. These processes are martingales with respect to the filtration
defined by the Monte Carlo stochasticity.
From the Doob--Kolmogorov martingale inequality [e.g., \citet{grimmett92}],
%
%
\begin{equation}
\prob\Bigl[
\sup_j |M^{n,k}_j|\ge\epsilon
\Bigr] \le
\frac{1}{\epsilon^2}\sum_{m=n}^{\infty}a_{m}^2
\sup_{|\theta|<k}\VarMC(D_m(\theta) ).
\label{eqkc2}
\end{equation}
Define events $F_n=\{\sup_j |M^{n}_j|\ge\epsilon\}$ and
$F_{n,k}=\{\sup_j |M^{n,k}_j|\ge\epsilon\}$.
It\vspace*{1.5pt} follows from \hyperlink{C4}{(B4)} and (\ref{eqkc2}) that
$\lim_{n\to\infty}\prob\{F_{n,k}\}=0$ for each $k$.
In light of the nondivergence assumed in \hyperlink{C2}{(B2)}, this
implies $\lim_{n\to\infty}\prob\{F_{n}\}=0$ which is
exactly (\ref{eqkc}).

To expand on this final assertion, let
$
\Omega=\{\sup_m |\hat\theta_m | < \infty\}
$
and
$
\Omega_k=\{\sup_m |\hat\theta_m | < k \}
$.
Assumption \hyperlink{C2}{(B2)} implies that $\prob(\Omega)=1$.
Since the sequence of events $\{\Omega_k\}$ is increasing up to
$\Omega$, we have $\lim_{k\rightarrow\infty}\prob(\Omega_k)=\prob
(\Omega)=1$.
Now observe that $\Omega_k\cap{F_{n,j}}=\Omega_k\cap{F_n}$ for all
$j\geq{k}$, as there is no truncation\vspace*{-1pt} of the sequence $\{\xi_m^j,
m=1,2,\ldots\}$ for outcomes in $\Omega_k$ when $j\geq{k}$.
Then,
\begin{eqnarray*}
\lim_{n\to\infty}\prob[F_n]&\le&\lim_{n\to\infty}\prob[F_n\cap
\Omega_k] + 1-\prob[\Omega_k]\\
&=&\lim_{n\to\infty}\prob[F_{n,k}\cap\Omega_k] + 1-\prob[\Omega
_k]\\
&\le&\lim_{n\to\infty}\prob[F_{n,k}] + 1- \prob[\Omega_k]\\
&=& 1-\prob[\Omega_k].
\end{eqnarray*}
Since $k$ can be chosen to make $1-\prob[\Omega_k]$ arbitrarily
small, it follows that $\lim_{n\to\infty}\prob[F_n]=0$.
\end{pf}

\subsection{Some standard results on sequential Monte Carlo and
importance sampling}
\label{appsmc}

A general convergence result on sequential Monte Carlo combining
results by \citet{crisan02} and \citet{delmoral01} is
stated in our
notation as Theorem~\ref{thcrisan} below.
The theorem is stated for a POMP model with generic density $\genf$,
parameter vector $\gentheta$, Markov process $\genX_{0\dvtx N}$\vspace*{1pt} taking
values in $\genRX^{N+1}$, observation process $\genY_{1\dvtx N}$ taking
values in $\genRY^{N}$, and data $\gendata_{1\dvtx N}$.
For application to the unperturbed model one sets\vspace*{1pt} $\genf=f$, $\genX
_n=X_n$, $\genRX=\RX$, $\genY_n=Y_n$ and $\gentheta
=\theta$.
For application to the stochastically perturbed model one sets $\genf
=g$, $\genX_n=(\chX_n,\chTheta_n)$, $\genRX=\RX\times\R
^\dimtheta$, $\genY_n=\chY_n$ and $\gentheta=(\theta
^{[N+1]},\sigma,\tau)$.
When applying Theorem~\ref{thcrisan} in the context of Theorem~\ref{th2}, the explicit expressions for the constants $\cA$ and $\cFF$
are required to show that the bounds in (\ref{eqphi1}) and (\ref{eqphibias}) apply uniformly for a collection of models indexed by
the approximation parameters $\{\tau_m\}$ and $\{\sigma_m\}$.



\begin{theorem}[{[\citet{crisan02};
\citet{delmoral01}]}] \label{thcrisan}
%
%
Let $\genf$ be a generic density for a POMP model having parameter
vector $\gentheta$, unobserved Markov process $\genX_{0\dvtx N}$,
observation process $\genY_{1\dvtx N}$ and data $\gendata_{1\dvtx N}$.
Define $\genX^\filter_{n,j}$ via applying Algorithm~\ref{algpf}
with $J$ particles.
Assume that $\genf_{\genY_n|\genX_n}(\gendata_n\given\tilde{x}
_n; \gentheta)$\vspace*{-2pt} is bounded as a function of $\tilde{x}_n$.
For any $\testFunction\dvtx \genRX\to\R$, denote the filtered mean of
$\testFunction(\genX_n)$ and its Monte Carlo estimate by
%
%
\begin{equation} \label{eqphi}
\testFunction^\filter_n=\int\testFunction(\tilde{x}_n)\genf
_{\genX
_n|\genY_{1\dvtx n}}(\tilde{x}_n\given\gendata_{1\dvtx n}; \gentheta
)\,
d\tilde{x}_n, \qquad
\widetilde\testFunction^\filter_n = \frac{1}{J}\sum_{j=1}^J
\testFunction(\genX^\filter_{n,j} ).
\end{equation}
There are constants $\cA$ and $\cFF$,
independent of $J$, such that
%
%
\begin{eqnarray}
\label{eqphi1}
\EMC[
(\widetilde\testFunction^\filter_n - \testFunction^\filter_n)^2
]
&\le&
\frac{\cA\sup_{\tilde{x}}|\testFunction(\tilde{x})|^2}{J},
\\
\label{eqphibias}
| \EMC[
\widetilde\testFunction^\filter_n - \testFunction^\filter_n
] |
&\le&
\frac{\cFF\sup_{\tilde{x}}|\testFunction(\tilde{x})|}{J}.
\end{eqnarray}
Specifically, $\cA$ and $\cFF$ can be written as linear functions of
$1$ and $\eta_{n,1},\ldots,\break\eta_{n,n}$ defined as
%
%
\begin{equation} \label{eqgamma}
\eta_{n,i} = \prod_{k=n-i+1}^{n} \biggl(
\frac{\sup_{\tilde{x}_k}\genf_{\genY_k|\genX_k}(\gendata
_k\given\tilde{x}_k; \gentheta)}{\genf_{\genY_k|\genY
_{1\dvtx k-1}}(\gendata_k\given\gendata_{1\dvtx k-1}; \gentheta)}
\biggr)^2.
\end{equation}
\end{theorem}

\begin{pf}
Theorem~2 of \citet{crisan02} derived (\ref{eqphi1}), and here we
start by focusing on the assertion that the constant $\cA$ in equation
(\ref{eqphi1}) can be written as a linear function of $1$ and the
quantities $\eta_{n,1},\ldots,\eta_{n,n}$ defined in (\ref{eqgamma}). This was not explicitly mentioned by \citet
{crisan02} but
is a direct consequence of their argument.
\citeauthor{crisan02} [(\citeyear{crisan02}), Section~V] constructed
the following recursion, for
which $c_{n|n}$ is the constant $\cA$ in equation (\ref{eqphi1}).
For $n=1,\ldots,N$ and $c_{0|0}=0$, define
%
%
\begin{eqnarray}\label{eqB1}
c_{n|n}&=& \bigl(
\sqrt{C}+\sqrt{{\varc_{n|n}}}
\bigr)^2,
\\\label{eqB2}
{\varc_{n|n}}&=&4c_{n|n-1} \biggl(\frac{\|\genf_{\genY
_n|\genX_n}\|}{\genf_{\genY_n|\genY_{1\dvtx n-1}}(\gendata
_n \given\gendata_{1\dvtx n-1}; \gentheta)} \biggr)^2,
\\
\label{eqB3}
c_{n|n-1}&=& \bigl(1+\sqrt{c_{n-1|n-1}} \bigr)^2,
\end{eqnarray}
where $\|\genf_{\genY_n|\genX_n}\|=\sup_{\tilde{x}_n}\genf
_{\genY_n|\genX_n}(\gendata_n \given\tilde{x}_n;
\gentheta)$.
Here, $C$ is a constant that depends on the resampling procedure but
not on the number of particles $J$.
Now, (\ref{eqB1})--(\ref{eqB3}) can be reformulated by routine
algebra as
%
%
\begin{eqnarray}\label{eqB4}
c_{n|n}&\leq&K_1+K_2 {\varc_{n|n}},\\\label{eqB5}
{\varc_{n|n}}&\leq&K_3 q_n c_{n|n-1},
\\\label{eqB6}
c_{n|n-1}&\leq&K_4+K_5 c_{n-1|n-1},
\end{eqnarray}
where
$
q_n=\|\genf_{\genY_n|\genX_n}\|^2 [ \genf_{\genY
_n|\genY_{1\dvtx n-1}}(\gendata_n \given\gendata_{1\dvtx n-1};
\gentheta) ]^{-2}
$
and $K_1,\ldots,K_5$ are constants which do not depend on $\genf$,
$\gentheta$, $\gendata_{1\dvtx N}$ or $J$.
Putting (\ref{eqB5}) and (\ref{eqB6}) into~(\ref{eqB4}),
%
%
\begin{eqnarray}\label{eqB7}
c_{n|n}&\leq&K_1 + K_2K_3q_n c_{n|n-1} \nonumber
\\[-8pt]
\\[-8pt]
&\le&K_1+K_2K_3K_4q_n + K_2K_3K_5 q_n c_{n-1|n-1}.
\nonumber
\end{eqnarray}
Since $\eta_{n,i}=q_n\eta_{n-1,i}$ for $i<n$, and $\eta_{n,n}=q_n$,
the required assertion follows from (\ref{eqB7}).\vadjust{\goodbreak}

To show (\ref{eqphibias}), we introduce the unnormalized filtered
mean $\testFunction^U_n$ and its Monte Carlo estimate $\widetilde
\testFunction^U_n$, defined by
%
%
\begin{equation}
\label{eqU1}
\testFunction^U_n= \testFunction^\filter_n \prod_{k=1}^n \genf
_{Y_k|Y_{1\dvtx k-1}}(\gendata_k\given\gendata_{1\dvtx k-1};
\gentheta),
\widetilde\testFunction^U_n=\widetilde\testFunction^\filter_n
\prod_{k=1}^n \frac{1}{J} \sum_{j=1}^J w(k,j),
\end{equation}
where $w(k,j)$ is computed in step~\hyperlink{step:w}{4} of
Algorithm~\ref{algpf} when evaluating $\widetilde\testFunction^\filter_n$.
Then, \citet{delmoral01} showed that
%
%
\begin{eqnarray}\label{eqgamma1} \label{equnbiased}
 \quad \EMC[\widetilde\testFunction^U_n ]
& = &
\testFunction^U_n,
\\\label{eqgamma2}
 \quad \EMC
[
(
\widetilde\testFunction^U_n-\testFunction^U_n
)^2
]
& \le&
\frac{({n+1})\sup_{\tilde{x}}|\testFunction(\tilde{x})|^2}{J}
\prod_{k=1}^n
\Bigl(
\sup_{\tilde{x}_k} \genf_{\genY_k|\genX_k}(\gendata
_k\given
\tilde{x}_k; \gentheta)
\Bigr)^2
.
\end{eqnarray}
We now follow an approach of \citeauthor{delmoral01} [(\citeyear{delmoral01}), equation~3.3.14],\vspace*{1pt} by
defining the unit function $1(\tilde{x})=1$ and observing that \rule
{0mm}{4mm} $\testFunction^F_n=\testFunction^U_n /1^U_n$ and
$\widetilde\testFunction^F_n=\widetilde\testFunction^U_n /
\widetilde1^U_n$. Then (\ref{equnbiased}) implies the identity
%
%
\begin{equation}\label{eqbias1}
\EMC[ \widetilde\testFunction^F_n -\testFunction^F_n ] =
\EMC\biggl[
(\widetilde\testFunction^F_n -\testFunction^F_n )
\biggl(1-\frac{\widetilde1^U_n}{1^U_n} \biggr)
\biggr].
\end{equation}
%
Applying the Cauchy--Schwarz inequality to (\ref{eqbias1}), making use
of (\ref{eqphi1}) and (\ref{eqgamma2}),
gives (\ref{eqphibias}).
\end{pf}

We now give a corollary to Theorem~\ref{thcrisan} for a latent
variable model $(\genX,\genY)$, as defined in Section~\ref{seciis},
having generic density $\genf$, parameter vector $\gentheta$,\vspace*{1pt}
unobserved variable $\genX$ taking values in $\genRX$, observed
variable $\genY$ taking values in $\genRY$, and data~$\gendata$.
Importance sampling for such a model is a special case of sequential
Monte Carlo, with $N=1$ and no resampling step.
We present and prove a separate result, which takes advantage of the
simplified situation, to make Section~\ref{seciis} and the proof of
Theorem~\ref{thsaiis} self-contained.
In the context of Theorem~\ref{thsaiis}, one sets $\genf=g$, $\genX
=(\chX,\chTheta)$, $\genRX=\RX\times\R^\dimtheta$, $\genY
_n=\chY$ and $\gentheta=(\theta,\tau)$.\looseness=1

\begin{corollary} \label{thcrisancor}
%
%
Let $\genf$ be a generic density for the latent variable model $(\genX
,\genY)$ with parameter vector $\gentheta$ and data $\data$.
Let $\{\genX_j, j=1,\ldots,J\}$ be $J$ independent Monte Carlo draws
from $\genf_{\genX}(\tilde{x}; \gentheta)$ and let $w_j=\genf
_{\genY
|\genX}(\data\given\genX_j; \gentheta)$.
Letting $\testFunction\dvtx \genRX\to\R$ be a bounded function, write
the conditional expectation of $\testFunction(\genX)$ and its
importance sampling estimate as
%
%
\begin{equation} \label{eqphiiis}
\testFunction^C = \int\testFunction(\tilde{x}_n)\genf_{\genX
|
\genY}(\tilde{x}\given\gendata; \gentheta)\,d\tilde{x},
\widetilde\testFunction^C = \frac{\sum_{j=1}^J w_j\testFunction
(\genX_j )}{\sum_{j=1}^J w_j}.
\end{equation}
Assume that $\genf_{\genY|\genX}(\gendata\given\tilde
{x};
\gentheta)$ is bounded as a function of $\tilde{x}$. Then,
%
%
\begin{eqnarray}
\label{eqphi1iis}
\VarMC(
\widetilde\testFunction^C
)
&\le&
\frac{4 \sup_{\tilde{x}}|\testFunction(\tilde{x})|^2
\sup_{\tilde{x}} (
\genf_{\genY|\genX}(\data\given\tilde{x}; \gentheta)
)^2
}{J (
\genf_{\genY}(\data; \gentheta)
)^2}
,
\\
\label{eqphibiasiis}
| \EMC[
\widetilde\testFunction^C ] - \testFunction^C
|
&\le&
\frac{2 \sup_{\tilde{x}}|\testFunction(\tilde{x})|
\sup_{\tilde{x}} (
\genf_{\genY|\genX}(\data\given\tilde{x}; \gentheta)
)^2
}{
J (
\genf_{\genY}(\data; \gentheta)
)^2}.
\end{eqnarray}
\end{corollary}
\begin{pf}
We introduce normalized weights $\hat w_j=w_j / \genf_{\genY
}(\data; \gentheta)$ and a normalized importance sampling estimator
$\hat\testFunction^C=\frac{1}{J}\sum_{j=1}^J \hat w_j\testFunction
(\genX_j )$\vspace*{-1pt} to compare to the self-normalized estimator in
(\ref{eqphiiis}).
It is not hard to check that\vspace*{1pt} $\EMC[\hat\testFunction^C
]=\testFunction^C$ and~$
\VarMC(\hat\testFunction^C )
\le
\frac{1}{J} (
\sup_{\tilde{x}} |\testFunction(\tilde{x})|
\sup_{\tilde{x}} \genf_{\genY|\genX}(\data\given\tilde
{x};
\gentheta)
/
\genf_{\genY}(\data; \gentheta)
)^2.
$
Now,\break
$
\VarMC(\widetilde\testFunction^C )\le2 \{\VarMC
(\widetilde\testFunction^C-\hat\testFunction^C )+\VarMC
(\hat\testFunction^C ) \}
$
and
\begin{eqnarray*}
\VarMC(\widetilde\testFunction^C-\hat\testFunction^C )
&\le&
\EMC\Biggl[ \Biggl( \frac{\fracd{1}{J}\sum_{j=1}^J \hat w_j
\testFunction(\genX_j)}{\fracd{1}{J}\sum_{j=1}^J \hat w_j} - \frac
{1}{J}\sum_{j=1}^J \hat w_j \testFunction(\genX_j) \Biggr)^{
2} \Biggr]
\\
&=& \EMC\Biggl[
\biggl(
\frac{\fracd{1}{J}\sum_{j=1}^J \hat w_j \testFunction(\genX
_j)}{\fracd{1}{J}\sum_{j=1}^J \hat w_j}
\biggr)^{ 2}
\Biggl(
1-\frac{1}{J}\sum_{j=1}^J \hat w_j
\Biggr)^{ 2}
\Biggr]
\\
&\le&
\sup_{\tilde{x}} |\testFunction(\tilde{x})|^2 \EMC
\Biggl[
1-\frac{1}{J}\sum_{j=1}^J \hat w_j
\Biggr]^2\\
&\le&
\frac{\sup_{\tilde{x}}|\testFunction(\tilde{x})|^2
\sup_{\tilde{x}}
(
\genf_{\genY|\genX}(\data\given\tilde{x}; \gentheta)
)^2
}{
J (
\genf_{\genY}(\data; \gentheta)
)^2}.
\end{eqnarray*}
This demonstrates (\ref{eqphi1iis}). To show (\ref{eqphibiasiis}), we write
\begin{eqnarray*}
|
\EMC[\widetilde\phi^C ]-\phi^C |
&=&
| \EMC[\widetilde\phi^C-\hat\phi^C ]
|\\
  &=&
\Biggl|
\EMC\Biggl[
(
\widetilde\phi^C - \EMC[ \widetilde\phi^C]
)
\Biggl(
1-\frac{1}{J}\sum_{j=1}^J \hat w_j
\Biggr)
\Biggr]
\Biggr|
\\
&\le&
\sqrt{ \VarMC(\widetilde\phi^C) \VarMC
\Biggl(
1-\frac{1}{J}\sum_{j=1}^J \hat w_j
\Biggr)
}
\\
  &\le&
\frac{{2} \sup_{\tilde{x}}|\testFunction(\tilde{x})|
\sup_{\tilde{x}}
(
\genf_{\genY|\genX}(\data\given\tilde{x}; \gentheta)
)^2
}{
J (
\genf_{\genY}(\data; \gentheta)
)^2}.
\end{eqnarray*}
\upqed
\end{pf}
\end{appendix}

\section*{Acknowledgments}
The authors are grateful for constructive advice from three anonymous
referees and the Associate Editor.


%

\printaddresses

\end{document}